\documentclass[12pt,reqno]{amsart}

\title[]{On Electroconvection in Porous Media}

\author{Elie Abdo}
\address{Department of Mathematics, Temple University, Philadelphia, PA 19122}
\email{abdo@temple.edu}
\author{Mihaela Ignatova}
\address{Department of Mathematics, Temple University, Philadelphia, PA 19122}
\email{ignatova@temple.edu}

\usepackage[margin=1in]{geometry}
\usepackage{amsmath, amsthm, amssymb}
\usepackage{times}
\usepackage{color}
\usepackage{hyperref}
\newcommand{\pa}{\partial}
\newcommand{\la}{\label}
\newcommand{\fr}{\frac}
\newcommand{\na}{\nabla}
\newcommand{\be}{\begin{equation}}
\newcommand{\ee}{\end{equation}}
\newcommand{\ba}{\begin{array}{l}}
\newcommand{\ea}{\end{array}}
\newcommand{\Rr}{{\mathbb R}}
\newtheorem{thm}{Theorem}
\newtheorem{prop}{Proposition}

\newtheorem{defi}{Definition}
\newcommand{\beg}{\begin}

\renewcommand{\l}{\Lambda}

\newtheorem{lem}{Lemma}

\usepackage{MnSymbol,wasysym}

\newcommand{\N}{\mathbb N}

\newcommand{\R}{\mathbb R}

\def\ZZ{{\mathbb Z}}

\def\PP{\mathbb P}

\date{today}
\begin{document}
\begin{abstract} We address existence, uniqueness and analyticity of solutions of an electroconvection model in porous media. 
\end{abstract} 

\vspace{.5cm}
\noindent

\noindent

\maketitle
\section{Introduction}\la{intro}
Electroconvection models describe the evolution of charge distributions in fluids \cite{turb}. In the case studied in \cite{exp,num}, the fluid occupies a thin region modeled as a two dimensional domain. The charge distribution is carried by the fluid and diffuses due to the parallel component of the electrical field.  This results in a nonlocal transport equation for the charge density $\rho$,
\be
\pa_t \rho + u \cdot \na \rho + \Lambda \rho = 0  
\la{1}
\ee 
where $\Lambda = (-\Delta)^{\fr{1}{2}}$ is the square root of the two dimensional Laplacian. The fluid is incompressible and is forced by electrical forces
\be
F = \rho E
\la{F}
\ee
where $E$ is the parallel component of the electrical field,
\be
E = -\na \Phi,
\la{E}
\ee
with $\na$ the gradient in $\Rr^2$. The relationship between the electrical potential $\Phi$ and the charge distribution confined to a two dimensional region is 
\be
\Phi = \Lambda^{-1}\rho
\la{phirho}
\ee
and we thus have
\be
F = -\rho R\rho
\la{frho}
\ee
with $R = \nabla \Lambda^{-1}$ the Riesz transforms. 
In general, the fluid obeys Navier-Stokes or related equations driven by the forces $F$. In this paper we consider flow through a porous medium, in which the dominant dissipation mechanism is due not to the viscosity of the fluid, but rather to an effective damping caused by flow through pores. The Stokes operator is then replaced by $u + \nabla p$.  We consider a system in which the fluid equilibrates rapidly and the Reynolds number is low, so that forces are balanced by damping, 
\be
u + \na p = F.
\la{uF}
\ee
This balance, together with \eqref{frho} and the requirement of incompressibility,
\be
\na\cdot u = 0,
\la{inco}
\ee
leads to
\be \la{2}
u = - \PP (\rho R \rho)
\ee
where $\mathbb{P}$ is the Leray-Hodge projector on divergence-free vector fields. In \cite{IS}, the balance law 
\eqref{2} was used to describe the solvent in a Nernst-Planck-Darcy system of ionic diffusion in 2D and 3D. 
The electroconvection  situation described above leads to the active scalar equation \eqref{1} with constitutive law \eqref{2}, which is the equation we study in this work.
The equation is $L^{\infty}$-critical, and resembles critical SQG \cite{CI2,CI3,CTV,I} except for the constitutive law \eqref{2} which in this case is nonlinear and doubly nonlocal. The well known global regularity of critical SQG \cite{caf,knv} is not available here.  In this paper we show that the equation \eqref{1}, \eqref{2} has global weak solutions. We describe local existence and uniqueness results for strong solutions. We also show that solutions with small initial data in Besov spaces slightly smaller than $L^{\infty}$ exist globally and are Gevrey regular. 

This paper is organized as follows. In section~\ref{Prelim}, we recall results about Besov spaces and  Littlewood-Paley decomposition. In section~\ref{weak}, we prove existence of global in time weak solutions of \eqref{1}, \eqref{2} for initial data in $L^{2+\delta}(\R^2)$ for some $\delta > 0$. If the initial data is in $L^p(\R^2)$ for $p \in (2, \infty]$, then the $L^p$ norm of any solution of  \eqref{1}, \eqref{2} remains bounded in time. Section~\ref{loc} is devoted to the existence of a unique local strong solution under higher regularity assumptions for the initial data. In section~\ref{BesovExistence}, we show that a global in time solution in Besov spaces exists provided that the initial data is sufficiently small in Besov spaces that are slightly smaller than $L^{\infty}(\R^2)$. In section~\ref{BesovAnalyt} we prove that solutions are Gevrey regular under a smallness condition imposed on the initial data. In section~\ref{reg}, we study the regularity and long time behavior of solutions for small initial data whereas in section~\ref{regarb}, we show that H\"{o}lder continuity of the charge distribution is a sufficient condition for the smoothness of solutions for arbitrary initial data, a result that is similar to the situation for SQG. In section~\ref{per}, we treat the periodic case, and we prove that the solution of the problem \eqref{1}, \eqref{2} posed on the two dimensional torus converges exponentially in time to zero. Finally, we consider in section~\ref{sub} the subcritical Darcy's law electroconvection, and we show existence of global smooth solutions for arbitrary initial data. 

\section{Preliminaries} \la{Prelim}

For $f \in \mathcal{S}'(\R^2)$, we denote the Fourier transform of $f$
by
\be 
\mathcal{F} f (\xi) = \widehat{f}(\xi) = \fr{1}{2\pi} \int_{\R^2} f(x) e^{-i \xi \cdot x} dx 
\ee
and its inverse by $\mathcal{F}^{-1}$.

Let $\Phi$ be a nonnegative, nonincreasing, infinitely differentiable, radial function such that $\Phi(r) = 1$ for $r \in \left[0, \fr{1}{2} \right]$ and $\Phi(r) = 0$ for $r \in \left[\fr{5}{8}, \infty \right]$. Let 
\be \la{dyadicpart} 
\Psi(r) = \Phi \left(\fr{r}{2} \right) - \Phi(r).
\ee 
For each $j \in \mathbb{Z}$, let 
\be \la{dyadicpart2}
\Psi_j (r) = \Psi (2^{-j}r).
\ee 
We have 
\be 
\Phi(|\xi|) + \sum\limits_{j=0}^{\infty} \Psi_j (|\xi|) = 1
\ee for all $\xi \in \R^2$ and 
\be 
\sum\limits_{j=-\infty}^{\infty} \Psi_j (|\xi|) = 1
\ee for all $\xi \in \R^2 \setminus \left\{0\right\}$. We define the homogeneous dyadic blocks 
\be 
\Delta_j f(x) = \fr{1}{2\pi} \int_{\R^2} \Psi_j (|\xi|) \widehat{f}(\xi) e^{i\xi \cdot x} d\xi 
= \mathcal{F}^{-1} \left[\Psi_j (|\cdot|) \widehat{f}(\cdot) \right](x)
\ee and the lower frequency cutoff functions 
\be 
S_j f = \sum\limits_{k \le j -1} \Delta_k f.
\ee
We note that the Fourier transform of each dyadic block is compactly supported. More precisely, we have 
\be \la{supports}
\mathrm{supp} \; \mathcal{F} (\Delta_j f) \subset 2^j \left[\fr{1}{2}, \fr{5}{4} \right]
\ee 
for all $j \in \mathbb{Z}$. 

Let $\mathcal{S}_h' (\mathbb{R}^2)$ be the set of all tempered distributions $u \in \mathcal{S}'(\R^2)$ such that 
\be 
\lim\limits_{j \to -\infty} S_j u = 0
\ee in $\mathcal{S}'(\R^2)$.
For $f \in \mathcal{S}_h'(\R^2)$, we denote the homogeneous Littlewood-Paley decomposition of $f$ by 
\be 
f = \sum\limits_{j \in \mathbb{Z}} \Delta_j f.
\ee

For $s \in \R, 1 \le p, q \le \infty$, we denote the homogeneous Besov space 
\be 
\dot{B}_{p,q}^{s}(\R^2) = \left\{ f \in \mathcal{S}_h'(\R^2) : \|f\|_{\dot{B}_{p,q}^s(\R^2)} < \infty \right\}
\ee where 
\be 
\|f\|_{\dot{B}_{p,q}^s(\R^2)} = \left(\sum\limits_{j \in \mathbb{Z}} 2^{jsq} \|\Delta_j f\|_{L^p(\mathbb{R}^2)}^q \right)^{1/q}
\ee 
and the inhomogeneous Besov space 
\be 
B_{p,q}^{s}(\R^2) = \left\{ f \in \mathcal{S}'(\R^2) : \|f\|_{{B}_{p,q}^s(\R^2)} < \infty \right\}
\ee where 
\be 
\|f\|_{{B}_{p,q}^s(\R^2)} = \left(2^{-sq} \|\tilde{\Delta}_{-1}f\|_{L^p(\mathbb{R}^2)}^q + \sum\limits_{j = 0}^{\infty} 2^{jsq} \|\Delta_j f\|_{L^p(\mathbb{R}^2)}^q \right)^{1/q}
\ee
with the usual modification when $q = \infty$. 
Here 
\be 
\tilde{\Delta}_{-1}f = \fr{1}{2\pi} \int_{\R^2} \Phi (|\xi|) \widehat{f}(\xi) e^{i\xi \cdot x} d\xi 
= \mathcal{F}^{-1} \left[\Phi (|\cdot|) \widehat{f}(\cdot) \right](x).
\ee
We note that the definition of the space $\dot{B}_{p,q}^{s}$ is independent of the function $\Phi$ in terms of which the dyadic blocks are defined. Indeed, any other dyadic partition yields an equivalent norm. 

If $s>0$, $1 \le p, q \le \infty$, then 
\be \la{hominhom}
B_{p,q}^s(\R^2) = \dot{B}_{p,q}^s(\R^2) \cap L^p(\R^2).
\ee Moreover, the norms $\|f\|_{{B}_{p,q}^s(\R^2)}$ and $\|f\|_{\dot{B}_{p,q}^s(\R^2)} + \|f\|_{L^p(\R^2)}$ are equivalent. 

We also consider the following time dependent homogeneous Besov spaces 
\be 
L^r (0, T; \dot{B}_{p,q}^s(\R^2)) = \left\{f (t)\in \mathcal{S}_h'(\R^2) : \|f\|_{L^r(0,T; \dot{B}_{p,q}^s(\R^2))} = \| \|f(\cdot, t)\|_{\dot{B}_{p,q}^s(\R^2)} \|_{L^r(0,T)} < \infty \right\}
\ee and 
\be 
\tilde{L}^r (0, T; \dot{B}_{p,q}^s(\R^2)) 
= \left\{f(t) \in \mathcal{S}_h'(\R^2) :  \|f\|_{\tilde{L}^r(0,T;\dot{B}_{p,q}^s(\R^2) )} < \infty \right\},
\ee
where
\[ \|f\|_{\tilde{L}^r(0,T;\dot{B}_{p,q}^s(\R^2) )} = \left(\sum\limits_{j \in \mathbb{Z}} 2^{jsq} \|\Delta_j f\|_{L^r(0, T; L^p(\mathbb{R}^2))}^q \right)^{1/q}.\]

We recall inequalities that are used in the paper (see for instance \cite{BCD,HK,T}).

\beg{prop} \la{mainprop} Let $f \in \mathcal{S}_h'(\R^2)$.
\begin{enumerate}
\item
(Bernstein's inequality) Let $1 \le p \le q \le \infty$. Let $k$ be a nonnegative integer. Then 
\be \la{BER1}
\sup\limits_{|\alpha|= k} \|\partial^{\alpha} \Delta_j f\|_{L^p(\mathbb{R}^2)} \le C_k 2^{jk} \|\Delta_j f\|_{L^p(\mathbb{R}^2)}
\ee holds for all $j \in \mathbb{Z}$.

\item 
Let $1 \le p \le q \le \infty$. Then \be \la{BER2} 
\|\Delta_j f\|_{L^q(\mathbb{R}^2)} \le C2^{2j \left(\fr{1}{p} - \fr{1}{q} \right)} \|\Delta_j f\|_{L^p(\mathbb{R}^2)}
\ee 
holds for all $j \in \mathbb{Z}$. Moreover, the continuous Besov embedding 
\be \la{BESOVEMBED}
\dot{B}_{p_1, q_1}^{s}(\mathbb{R}^2) \hookrightarrow \dot{B}_{p_2, q_2}^{s - 2\left(\fr{1}{p_1} - \fr{1}{p_2} \right)}(\mathbb{R}^2)
\ee  holds for $1 \le p_1 \le p_2 \le \infty, 1 \le q_1 \le q_2 \le \infty$ and $s \in \R$.

\item 
Let $1 \le p \le \infty, t \ge 0, \alpha > 0$.  Then 
\be  \la{localization}
\|e^{-t \Lambda^{\alpha}} \Delta_j f\|_{L^p(\mathbb{R}^2)} \le Ce^{-C^{-1}t2^{j\alpha}} \|\Delta_j f\|_{L^p(\mathbb{R}^2)}
\ee holds for all $j \in \mathbb{Z}$. Here $\Lambda^{\alpha}$ is the fractional Laplacian of order $\alpha$ defined as a Fourier multiplier with symbol $|\xi|^{\alpha}$.

\item 
Let $R = (R_1, R_2)$ be the Riesz transform, i.e., for $k \in \left\{1, 2\right\}$,  $R_k = \pa_k \l^{-1}$. For each $p \in [1, \infty]$, there is a positive constant $C > 0$ depending only on $p$ (independent of $j$) such that 
\be \la{jRB}
\|\Delta_j R f\|_{L^p(\mathbb{R}^2)} \le C \|\Delta_j f\|_{L^p(\mathbb{R}^2)}
\ee holds for all  $j \in \mathbb{Z}$. Hence, for $s \in \R$ and $1\le p, q \le \infty$, $R$ is bounded from $\dot{B}_{p,q}^{s} (\mathbb{R}^2)$ to itself.  
\end{enumerate}
\end{prop}

The following decomposition formula holds.

\beg{prop} \la{Bon}
Let $f, g \in \mathcal{S}_h'(\R^2)$. Then 
\beg{align} \la{bonydecomp}
\Delta_j (fg) &= \sum\limits_{k \ge j-2} \Delta_j (S_{k+1}f \Delta_k g) 
+  \sum\limits_{k \ge j-2} \Delta_j (S_{k}g \Delta_k f)\nonumber
\\&=  \sum\limits_{k \ge j-2} \Delta_j (S_{k+1}g \Delta_k f) 
+  \sum\limits_{k \ge j-2} \Delta_j (S_{k}f \Delta_k g)
\end{align} holds for any $j \in \ZZ$. 
\end{prop}

The proof is based on Bony's paraproduct, and is presented in Appendix A.

Throughout this paper $C$ (or $C_i, i=1,2,\dots$) denotes a positive constant  that may change from line to line in the proofs. 

\section{Existence of Global Weak Solutions} \la{weak}
We consider the transport and nonlocal diffusion equation  
\be 
\pa_t \rho + u \cdot \na \rho + \Lambda \rho = 0
\la{pde}
\ee
in the whole space $\R^2$, where
\be 
u = - \PP (\rho R \rho).
\la{udef}
\ee 
The initial data are
\be 
\rho(x,0) = \rho_0(x).
\la{initial}
\ee
Here $\mathbb{P}$ is the Leray-Hodge projector, $\Lambda = (-\Delta)^{\fr{1}{2}}$ is the fractional Laplacian, and $R = \nabla \Lambda^{-1}$ is the 2D vector of Riesz transforms.  

\begin{defi}
A solution $\rho$ of the initial value problem \eqref{pde}--\eqref{initial} is said to be a weak solution on $[0,T]$ if 
\be 
\rho \in L^{\infty}(0,T;L^2(\R^2)) \cap L^2(0,T;\dot{H}^{\fr{1}{2}}(\R^2))
\ee and $\rho$ obeys  
\be 
(\rho (t), \Phi)_{L^2} - (\rho_0, \Phi)_{L^2} 
- \int_{0}^{t} (\rho, u\cdot \na \Phi)_{L^2}  ds 
+ \int_{0}^{t} (\Lambda^{\fr{1}{2}} \rho, \Lambda^{\fr{1}{2}}  \Phi)_{L^2} ds 
= 0
\ee for all $\Phi \in H^{\fr{5}{2}}(\R^2)$ and a.e. $t \in [0,T]$.
\end{defi}

For $\epsilon \in (0,1]$, let $J_{\epsilon}$ be the standard mollifier operator $J_{\epsilon}f = J_{\epsilon} * f$, and let $\rho^{\epsilon}$ be the solution of 
\be \la{approx}
\pa_t \rho^{\epsilon} + \widetilde{u}^{\epsilon} \cdot \na \rho^{\epsilon} + \l \rho^{\epsilon} - \epsilon \Delta \rho^{\epsilon} = 0
\ee where
\be 
\widetilde{u}^{\epsilon} = - J_{\epsilon} \mathbb{P}(\rho^{\epsilon} R \rho^{\epsilon})
\la{udefapprox}
\ee
with smoothed out initial data
\be 
\rho_0^{\epsilon} = J_{\epsilon} \rho_0
\la{initialapprox}
\ee

\beg{rem}
We note that $\mathbb{P}$ and $J_{\epsilon}$ commutes, hence $\widetilde{u}^{\epsilon}$ is divergence free. 
\end{rem}

\begin{thm} \la{Weak1}
Let $T > 0$ be arbitrary. Let $\rho_0 \in L^2(\R^2)$. Then for each $\epsilon  \in (0,1]$, the mollified initial value problem \eqref{approx}--\eqref{initialapprox} has a solution $\rho^{\epsilon}$ on $[0,T]$ satisfying 
\be
\fr{1}{2} \|\rho^{\epsilon}(t)\|_{L^2}^2 + \int_{0}^{t} \|\Lambda^{\fr{1}{2}} \rho^{\epsilon} (s)\|_{L^2}^2 ds   \leq \fr{1}{2} \|\rho_0\|_{L^2}^2
\la{reg1}
\ee 
for all $t \in [0,T]$.
Moreover, the sequence $\left\{\rho^{1/n} \right\}_{n=1}^{\infty}$ has a subsequence that converges strongly in $L^2(0,T; L^2(\R^2))$ and weakly in $L^2(0,T; H^{\fr{1}{2}}(\R^2))$ to a function $\rho$ obeying 
\be  
\fr{1}{2} \|\rho (t)\|_{L^2}^2 + \int_{0}^{t} \|\Lambda^{\fr{1}{2}} \rho (s)\|_{L^2}^2 ds   \leq \fr{1}{2} \|\rho_0\|_{L^2}^2
\la{reg2}
\ee 
for a.e. $t \in [0,T]$. If $\rho_0 \in L^{2+ \delta}(\R^2)$ for some $\delta > 0$, then $\rho$ is a weak solution of \eqref{pde}--\eqref{initial} on $[0,T]$.
\end{thm}
The proof is found in Appendix B. 

As a consequence of the C\'ordoba-C\'ordoba inequality \cite{CC}, the $L^p$ norm of any solution of the equation \eqref{pde}--\eqref{udef} is bounded by the $L^p$ norm of the initial data for any $p \in (2, \infty]$:

\begin{prop} \la{crit}
Let $p > 2$ and $\rho_0 \in L^p(\R^2)$. Suppose $\rho$ is a smooth solution of \eqref{pde}--\eqref{initial} on $[0,T]$. Then 
\be 
\|\rho(t)\|_{L^p} \leq \|\rho_0\|_{L^p}
\la{decayp}
\ee 
holds for all $t \in [0,T]$. Moreover, if $\rho_0 \in L^{\infty}(\R^2)$, then 
\be 
\|\rho(t)\|_{L^{\infty}} \leq \|\rho_0\|_{L^{\infty}} 
\la{decayinf}
\ee 
holds for all $t \in [0,T]$.
\end{prop}

\textbf{Proof:} We multiply \eqref{pde} by $\rho|\rho|^{p-2}$ and we integrate in the space variable. We obtain the differential inequality 
\be 
\fr{d}{dt} \|\rho\|_{L^p}  \leq 0.
\ee 
This gives \eqref{decayp}. 
Letting $p \rightarrow \infty$, we obtain \eqref{decayinf}.

\beg{rem} \la{ltime} If $\rho$ is a smooth solution of \eqref{pde}--\eqref{initial} on $[0,T)$ and $\rho (\cdot, t) \in H^s(\R^2)$ for some $s > 1$ and for a.e. $t \in [0,T)$, then 
\be 
\|\rho(.,t)\|_{L^{\infty}} \le \fr{\|\rho_0\|_{L^{\infty}}}{1 + Ct\|\rho_0\|_{L^{\infty}}}
\ee for $t \in [0,T)$ (see \cite{CC}). This bound is useful to study the long time behavior of solutions. 
\end{rem}

\section{Existence of Local Strong Solutions} \la{loc}
  
\begin{defi}
A weak solution $\rho$ of \eqref{pde}--\eqref{initial} is said to be a strong solution  on $[0,T]$ if it obeys 
\be 
\rho \in L^{\infty} (0,T; \dot{H}^2(\R^2)) \cap  L^2(0,T; \dot{H}^{\fr{5}{2}}(\R^2)).
\ee 
\end{defi}

\begin{thm}\la{strongloc}
Let $\rho_0 \in H^2(\R^2)$. Then there exists $T_0 > 0$ depending only on $\|\rho_0\|_{H^2}$ such that a unique strong solution of the equation \eqref{pde}--\eqref{udef} exists on $[0,T_0]$.
\end{thm}

The proof is found in Appendix~\ref{appstrong}.

\section{Existence of Global Solutions in Besov Spaces} \la{BesovExistence}

In this section, we show the existence of a global in time solution in Besov spaces for sufficiently small initial data. The proof uses methods of \cite{BBT,C}.

\beg{thm} \la{solinBesov} Let $1 \le p < \infty$. Let $\rho_0 \in \dot{B}_{p,1}^{\fr{2}{p}}(\R^2)$ be sufficiently small. We consider the functional space $E_p$ defined by 
\be 
\la{Epdef}
E_p = \left\{f(t) \in \mathcal{S}_h'(\R^2) : \|f\|_{E_p} = \|f\|_{\tilde{L}_t^{\infty} \dot{B}_{p,1}^{\fr{2}{p}}} + \|f\|_{\tilde{L}_t^{1} \dot{B}_{p,1}^{\fr{2}{p} +1 }} < \infty \right\}.
\ee 
Then the equation \eqref{pde}--\eqref{udef} has a unique global in time solution $\rho \in E_p$.
\end{thm}

\textbf{Proof:} Let $\rho^{(0)} = 0$. For each positive integer $n$, let $\rho^{(n)}$ be the solution of 
\be \la{approxeq}
\beg{cases} \pa_t \rho^{(n)} + \l \rho^{(n)} = - u^{(n-1)}\cdot \na \rho^{(n-1)}, 
\\ u^{(n-1)} = -\mathbb{P} (\rho^{(n-1)}R \rho^{(n-1)}),
\\ \rho_0^{(n)} = \rho^{(n)} (x, 0) = \rho_0 
\end{cases}
\ee
posed on $\mathbb{R}^2$. We write $\rho^{(n)}$ in the integral form, 
\beg{align} \la{approxsol}
\rho^{(n)}(t) &= e^{-t\l}\rho_0 - \int_{0}^{t} e^{-(t-s)\l} \na \cdot (u^{(n-1)} \rho^{(n-1)})(s) ds \nonumber
\\&= e^{-t\l}\rho_0 - \mathcal{B}(u^{n-1}, \rho^{n-1})
\end{align} 
where $\mathcal{B}$ is the bilinear form defined by 
\be 
\mathcal{B}(v, \theta) = \int_{0}^{t} e^{-(t-s)\l} \na \cdot (v \theta)(s) ds. 
\ee 

\textit{Step 1.} Fix a positive integer $n$.  We show that 
\be \la{REC}
\|\rho^{(n)}\|_{E_p} \le C_1 \|\rho_0\|_{\dot{B}_{p,1}^{\fr{2}{p}}} + C_2 \|\rho^{(n-1)}\|_{E_p}^3.
\ee
We start by estimating $e^{-t\l}\rho_0$ in $E_p$. We apply $\Delta_j$ and we take the $L^p$ norm. In view of the bound \eqref{localization}, we have 
\be \la{L1}
\|e^{-t\l} \Delta_j \rho_0 \|_{L^p} \le Ce^{-C^{-1}t 2^{j}} \|\Delta_j \rho_0 \|_{L^p},
\ee hence 
\be 
\|e^{-t\l} \rho_0\|_{E_p} = \|e^{-t\l} \rho_0\|_{\tilde{L}_t^{\infty} \dot{B}_{p,1}^{\fr{2}{p}}} + \|e^{-t\l} \rho_0\|_{\tilde{L}_t^{1} \dot{B}_{p,1}^{\fr{2}{p} + 1}} 
\le C\|\rho_0\|_{\dot{B}_{p,1}^{\fr{2}{p}}}.
\ee
Now, we estimate the term $\mathcal{B}(u^{(n-1)}, \rho^{(n-1)})$ in $E_p$. First, we note that 
\be \la{BUV}
\|\mathcal{B}(u^{(n-1)}, \rho^{(n-1)})\|_{E^p} \le C \|u^{(n-1)} \rho^{(n-1)}\|_{\tilde{L}_t^1 \dot{B}_{p,1}^{\fr{2}{p}+1} }.
\ee
Indeed, we apply $\Delta_j$ to $\mathcal{B}(u^{(n-1)}, \rho^{(n-1)})$ and we estimate. On one hand, 
\beg{align}
\|\Delta_j \mathcal{B}(u^{(n-1)}, \rho^{(n-1)})\|_{L_t^{\infty}L^p} 
&\le C2^j \left\|\int_{0}^{t} e^{-c^{-1}(t-s)2^{j}} \|\Delta_j(u^{(n-1)}\rho^{(n-1)})(s)\|_{L^p} ds  \right\|_{L_t^{\infty}} \nonumber
\\&\le C2^j \|\Delta_j (u^{(n-1)} \rho^{(n-1)})\|_{L_t^1 L^p} 
\end{align}
in view of Bernstein's inequality \eqref{BER1} and the bound \eqref{localization}. We multiply by $2^{j\fr{2}{p}}$ and we take the $\ell^1$ norm. We obtain the bound 
\be \la{BUV1}
\|\mathcal{B} (u^{(n-1)}, \rho^{(n-1)})\|_{\tilde{L}_t^{\infty} \dot{B}_{p,1}^{\fr{2}{p}}} \le C \|u^{(n-1)} \rho^{(n-1)}\|_{\tilde{L}_t^1 \dot{B}_{p,1}^{\fr{2}{p}+1} }. 
\ee 
On the other hand, 
\beg{align} 
&\|\Delta_j \mathcal{B} (u^{(n-1)}, \rho^{(n-1)})\|_{L_t^1 L^p}
\le C \left\|\int_{0}^{t} 2^j e^{-c^{-1}(t-s)2^{j}} \|\Delta_j(u^{(n-1)}\rho^{(n-1)})(s)\|_{L^p} ds  \right\|_{L_t^{1}} \nonumber
\\&\le C \int_{0}^{\infty} \left(\int_{0}^{\infty} 2^j e^{-c^{-1}(t-s)2^{j}} \chi_{[0,t]}(s) dt \right)   \|\Delta_j(u^{(n-1)}\rho^{(n-1)})(s)\|_{L^p} ds \nonumber
\\&\le C \|\Delta_j (u^{(n-1)} \rho^{(n-1)})\|_{L_t^1 L^p}
\end{align} where $\chi_{E}$ denotes the characteristic function of the set $E$. 
Multiplying by $2^{j \left(\fr{2}{p} + 1 \right)}$ and taking the $\ell^1$ norm yields the bound 
\be \la{BUV2}
\|\mathcal{B} (u^{(n-1)}, \rho^{(n-1)})\|_{\tilde{L}_t^{1} \dot{B}_{p,1}^{\fr{2}{p}}} \le C \|u^{(n-1)} \rho^{(n-1)}\|_{\tilde{L}_t^1 \dot{B}_{p,1}^{\fr{2}{p}+1} }. 
\ee 
Combining \eqref{BUV1} and \eqref{BUV2}, we obtain \eqref{BUV}. Accordingly, our next goal is to show that 
\be \la{REC1}
\|u^{(n-1)} \rho^{(n-1)}\|_{\tilde{L}_t^1 \dot{B}_{p,1}^{\fr{2}{p}+1} } \le C\|\rho^{(n-1)}\|_{E^p}^3
\ee which gives \eqref{REC}.
In order to establish the bound \eqref{REC1}, we use the decomposition \eqref{bonydecomp}  
\be 
\Delta_j (u^{(n-1)}\rho^{(n-1)}) 
= \sum\limits_{k \ge j -2} \Delta_j (S_ku^{(n-1)} \Delta_k \rho^{(n-1)}) 
+ \sum\limits_{k \ge j-2} \Delta_j (S_{k+1} \rho^{(n-1)} \Delta_k u^{(n-1)}).
\ee
We apply the $L_t^1 L^p$ norm, we use the bound
\be 
\|\Delta_j f\|_{L^p} \le C \|f\|_{L^p}
\ee that holds for any $f \in \mathcal{S}_h'$ where $C$ is a positive universal constant independent of $j$, and we obtain 
\beg{align} \la{L2}
\|\Delta_j (u^{(n-1)}\rho^{(n-1)})\|_{L_t^1 L^p} 
&\le C\sum\limits_{k \ge j -2} \|S_ku^{(n-1)}\|_{L_t^{\infty}L^{\infty}} \|\Delta_k \rho^{(n-1)}\|_{L_t^1 L^p} \nonumber
\\&+ C\sum\limits_{k \ge j-2} \|S_{k+1} \rho^{(n-1)}\|_{L_t^{\infty}L^{\infty}}\| \Delta_k u^{(n-1)}\|_{L_t^1 L^p}.
\end{align}
In view of Bernstein's inequality \eqref{BER2}, we have 
\be \la{B11} 
\|S_{k+1} \rho^{(n-1)}\|_{L_t^{\infty}L^{\infty}} 
\le \sum\limits_{l \le k} \|\Delta_l \rho^{(n-1)}\|_{L_t^{\infty}L^{\infty}}
\le C\sum\limits_{l \le k} 2^{l\fr{2}{p}} \|\Delta_l \rho^{(n-1)}\|_{L_t^{\infty}L^p}
\le C\|\rho^{(n-1)}\|_{\tilde{L}_t^{\infty} \dot{B}_{p,1}^{\fr{2}{p}}}.
\ee
We show below that 
\be \la{nonlinear1}
\|S_k u^{(n-1)}\|_{L_t^{\infty}L^{\infty}} \le C\|\rho^{(n-1)}\|_{\tilde{L}_t^{\infty} \dot{B}_{p,1}^{\fr{2}{p}}}^2
\ee 
and 
\be \la{nonlinear2}
\|\Delta_k u^{(n-1)}\|_{L_t^1 L^p} \le C\|\rho^{(n-1)}\|_{\tilde{L}_t^{\infty} \dot{B}_{p,1}^{\fr{2}{p}}}\left( \sum\limits_{m \ge k-2} \|\Delta_m \rho^{(n-1)}\|_{L_t^1 L^p}  \right).
\ee
Using the bounds \eqref{nonlinear1} and \eqref{nonlinear2}, we obtain
\be \la{final}
\|\Delta_j (u^{(n-1)}\rho^{(n-1)})\|_{L_t^1 L^p}  \le C \|\rho^{(n-1)}\|_{\tilde{L}_t^{\infty} \dot{B}_{p,1}^{\fr{2}{p}}}^2 \left\{ \sum\limits_{k \ge j-2} \|\Delta_k \rho^{(n-1)}\|_{L_t^1 L^p} +  \sum\limits_{k \ge j-2}\sum\limits_{m \ge k-2} \|\Delta_m \rho^{(n-1)}\|_{L_t^1 L^p}  \right\}
\ee
We multiply \eqref{final} by $2^{j \left(\fr{2}{p} +1 \right)}$ and we take the $\ell^1$ norm. In view of Young's convolution inequality, we have in the first term
\beg{align} \la{final1}
&\sum\limits_{j \in \mathbb{Z}} \sum\limits_{k \ge j-2} 2^{j \left(\fr{2}{p} + 1 \right)} \|\Delta_k \rho^{(n-1)}\|_{L_t^1 L^p}
= \sum\limits_{j \in \mathbb{Z}} \sum\limits_{k \ge j-2} 2^{-(k-j)  \left(\fr{2}{p} + 1 \right)} 2^{k \left(\fr{2}{p} + 1 \right)}\|\Delta_k \rho^{(n-1)}\|_{L_t^1 L^p} \nonumber
\\&\le \left(\sum\limits_{j \ge -2} 2^{-j \left(\fr{2}{p} +1 \right)} \right) \left(\sum\limits_{j \in \mathbb{Z}} 2^{j \left(\fr{2}{p} + 1 \right)} \|\Delta_j \rho^{(n-1)}\|_{L_t^1 L^p} \right)
\le C\|\rho^{(n-1)}\|_{\tilde{L}_t^1 \dot{B}_{p,1}^{\fr{2}{p}+1}}.
\end{align} 
For the second summation on the right hand side of \eqref{final}, we apply Fubini's theorem and then we estimate as  in \eqref{final1}. Thus, we have 
\beg{align} \la{final2}
&\sum\limits_{j \in \mathbb{Z}} \sum\limits_{k \ge j-2}  \sum\limits_{m \ge k-2} 2^{j \left(\fr{2}{p} + 1 \right)} \|\Delta_m \rho^{(n-1)}\|_{L_t^1 L^p} \nonumber
\\&= \sum\limits_{j \in \mathbb{Z}} \sum\limits_{m \ge j-4}  \sum\limits_{j-2 \le k \le m+2} 2^{-(m-j) \left(\fr{2}{p} + 1 \right)} 2^{m \left(\fr{2}{p} + 1 \right)}\|\Delta_m \rho^{(n-1)}\|_{L_t^1 L^p} \nonumber
\\&= \sum\limits_{j \in \mathbb{Z}} \sum\limits_{m \ge j-4}  (m-j+5) 2^{-(m-j) \left(\fr{2}{p} + 1 \right)} 2^{m \left(\fr{2}{p} + 1 \right)}\|\Delta_m \rho^{(n-1)}\|_{L_t^1 L^p} \nonumber
\\&\le C\sum\limits_{j \in \mathbb{Z}} \sum\limits_{m \ge j-4} 2^{-(m-j)\left(\fr{1}{p} + \fr{1}{2} \right)} 2^{m \left(\fr{2}{p} + 1 \right)}\|\Delta_m \rho^{(n-1)}\|_{L_t^1 L^p}  \nonumber
\\&+ 5 \sum\limits_{j \in \mathbb{Z}} \sum\limits_{m \ge j-4}  2^{-(m-j) \left(\fr{2}{p} + 1 \right)} 2^{m \left(\fr{2}{p} + 1 \right)}\|\Delta_m \rho^{(n-1)}\|_{L_t^1 L^p} \nonumber
\\&\le C\|\rho^{(n-1)}\|_{\tilde{L}_t^1 \dot{B}_{p,1}^{\fr{2}{p}+1}}.
\end{align} Here, we have used the fact that $x2^{-x} \le C2^{-\fr{x}{2}}$ for all $x \in \R$.
Putting \eqref{final1} and \eqref{final2} together, we obtain \eqref{REC1}. 

We end the proof of Step 1 by showing the estimates \eqref{nonlinear1} and \eqref{nonlinear2}. For each $l \in \mathbb{Z}$, we use again paraproducts to decompose $\Delta_l (\rho^{(n-1)}R\rho^{(n-1)})$ as 
\be \la{para}
\Delta_l (\rho^{(n-1)}R \rho^{(n-1)}) 
= \sum\limits_{m \ge l - 2} \Delta_l (S_{m+1} \rho^{(n-1)} \Delta_m R\rho^{(n-1)})
+ \sum\limits_{m \ge l - 2} \Delta_l (S_m R\rho^{(n-1)} \Delta_m \rho^{(n-1)}).
\ee
In view of the boundedness of the Riesz transform \eqref{jRB} and the definition of the Leray projector as 
\be 
\mathbb{P} = I + R \otimes R,
\ee we bound
\beg{align}
\|S_k u^{(n-1)}\|_{L_t^{\infty}L^{\infty}} 
&\le \sum\limits_{l \le k-1} \|\Delta_l u^{(n-1)}\|_{L_t^{\infty}L^{\infty}}
\le C\sum\limits_{l \le k-1} 2^{l \fr{2}{p}} \|\Delta_l u^{(n-1)}\|_{L_t^{\infty}L^p} \nonumber
\\&\le C\sum\limits_{l \le k-1}  2^{l \fr{2}{p}} \|\Delta_l (\rho^{(n-1)}R\rho^{(n-1)})\|_{L_t^{\infty}L^p}
\end{align} for any $p \in [1,\infty]$ and using the paraproduct decomposition \eqref{para}, we obtain 
\beg{align} 
\|S_k u^{(n-1)}\|_{L_t^{\infty}L^{\infty}}  
&\le C\sum\limits_{l \le k-1} 2^{l \fr{2}{p}} \sum\limits_{m \ge l - 2} \|S_{m+1} \rho^{(n-1)}\|_{L_t^{\infty} L^{\infty}} \|\Delta_m R\rho^{(n-1)}\|_{L_t^{\infty}L^p} \nonumber 
\\&+ C\sum\limits_{l \le k-1} 2^{l \fr{2}{p}} \sum\limits_{m \ge l - 2} \|S_m R\rho^{(n-1)}\|_{L_t^{\infty} L^{\infty}} \|\Delta_m \rho^{(n-1)}\|_{L_t^{\infty}L^p}.
\end{align} We note that
\be 
\|S_{m+1} \rho^{(n-1)}\|_{L_t^{\infty} L^{\infty}} \le C\|\rho^{(n-1)}\|_{\tilde{L}_t^{\infty} \dot{B}_{p,1}^{\fr{2}{p}}}
\ee 
as shown in \eqref{B11}. Moreover, in view of  \eqref{jRB}, we have
\beg{align} \la{N1}
\|S_m R\rho^{(n-1)}\|_{L_t^{\infty} L^{\infty}} 
&\le \sum\limits_{z \le m-1} \|\Delta_{z} R\rho^{(n-1)}\|_{{L_t^{\infty} L^{\infty}}}
\le C\sum\limits_{z \le m-1} 2^{z \fr{2}{p}} \|\Delta_{z} R\rho^{(n-1)}\|_{{L_t^{\infty} L^{p}}} \nonumber
\\&\le C\sum\limits_{z \le m-1} 2^{z \fr{2}{p}} \|\Delta_{z} \rho^{(n-1)}\|_{{L_t^{\infty} L^{p}}}
\le C\|\rho^{(n-1)}\|_{\tilde{L}_t^{\infty} \dot{B}_{p,1}^{\fr{2}{p}}}.
\end{align} Now we use the assumption that $p < \infty$ which implies that $\fr{2}{p} >0$ and so we can apply Young's convolution inequality to obtain 
\beg{align} \la{N2}
&\|S_k u^{(n-1)}\|_{L_t^{\infty}L^{\infty}}  
\le C\|\rho^{(n-1)}\|_{\tilde{L}_t^{\infty} \dot{B}_{p,1}^{\fr{2}{p}}} \left\{\sum\limits_{l \le k-1} 2^{l \fr{2}{p}} \sum\limits_{m \ge l-2} \|\Delta_m \rho^{(n-1)}\|_{L_t^{\infty}L^p}  \right\} \nonumber
\\&= C\|\rho^{(n-1)}\|_{\tilde{L}_t^{\infty} \dot{B}_{p,1}^{\fr{2}{p}}} \left\{\sum\limits_{l \le k-1} \sum\limits_{m \ge l-2}  2^{- (m-l) \fr{2}{p}} 2^{m\fr{2}{p}} \|\Delta_m \rho^{(n-1)}\|_{L_t^{\infty}L^p}  \right\}
\le C\|\rho^{(n-1)}\|_{\tilde{L}_t^{\infty} \dot{B}_{p,1}^{\fr{2}{p}}}^2 
\end{align} which proves \eqref{nonlinear1}. 
We proceed to show \eqref{nonlinear2}. Using the paraproduct decomposition \eqref{para} and the bound \eqref{jRB}, we have
\beg{align}
\|\Delta_k u^{(n-1)}\|_{L_t^1 L^p} 
&\le C\|\Delta_k (\rho^{(n-1)}R \rho^{(n-1)})\|_{L_t^1 L^p} \nonumber
\\&\le C \sum\limits_{m \ge k - 2} \|S_{m+1} \rho^{(n-1)}\|_{L_t^{\infty}L^{\infty}} \|\Delta_m R\rho^{(n-1)}\|_{L_t^1 L^p} \nonumber
\\&+ C\sum\limits_{m \ge k - 2} \|S_m R\rho^{(n-1)}\|_{L_t^{\infty} L^{\infty}} \|\Delta_m \rho^{(n-1)}\|_{L_t^1 L^p} \nonumber 
\\&\le C\|\rho^{(n-1)}\|_{\tilde{L}_t^{\infty} \dot{B}_{p,1}^{\fr{2}{p}}} \left(\sum\limits_{m \ge k-2} \|\Delta_m \rho^{(n-1)}\|_{L_t^1 L^p}  \right)
\end{align} yielding \eqref{nonlinear2}. This ends the proof of Step 1. 

\textit{Step 2.} We show that there exists an $\epsilon > 0$ sufficiently small such that if $C_1\|\rho_0\|_{\dot{B}_{p,1}^{\fr{2}{p}}} < \epsilon$, then the sequence $\left\{\rho^{(n)}\right\}_{n=1}^{\infty}$ converges to a unique solution $\rho$  of \eqref{pde}--\eqref{initial} obeying $\|\rho\|_{E_p} < 2\epsilon$. 

First, choose an $\epsilon > 0$ such that $C_2(2\epsilon)^3 < \epsilon$, where $C_2$ is the constant in \eqref{REC}, and suppose  that $C_1\|\rho_0\|_{\dot{B}_{p,1}^{\fr{2}{p}}} < \epsilon$.
Then an inductive argument yields 
\be  \la{binduc}
\|\rho^{(n)}\|_{E_p} < 2\epsilon 
\ee for all $n \ge 1$. Indeed, 
\be 
\|\rho^{(1)}\|_{E_p} \le C_1\|\rho_0\|_{\dot{B}_{p,1}^{\fr{2}{p}}} < \epsilon < 2\epsilon
\ee in view of \eqref{REC}. Suppose that 
\be 
\|\rho^{(n-1)}\|_{E_p} < 2\epsilon.
\ee Then 
\be 
\|\rho^{(n)}\|_{E_p} < \epsilon + C_2(2\epsilon)^3 < \epsilon + \epsilon = 2\epsilon.
\ee Therefore, we obtain \eqref{binduc}

Now, we show that the sequence  $\left\{\rho^{(n)}\right\}_{n=1}^{\infty}$ is Cauchy. Indeed, the difference $\rho^{(n)} - \rho^{(n-1)}$ obeys 
\beg{align} 
(\rho^{(n)} - \rho^{(n-1)}) (t) 
&= \int_{0}^{t} e^{-(t-s)\l} \na \cdot \left[u^{(n)}(\rho^{(n)} - \rho^{(n-1)}) - (u^{(n-1)} - u^{(n)})\rho^{(n-1)} \right](s) ds \nonumber
\\&= \mathcal{B}(u^{(n)}, \rho^{(n)} - \rho^{(n-1)}) - \mathcal{B}(u^{(n-1)} - u^{(n)}, \rho^{(n-1)}).
\end{align} As in Step 1 and using \eqref{binduc}, it can be shown that 
\beg{align} \la{cauchy}
\|\rho^{(n)} - \rho^{(n-1)} \|_{E_p} 
&\le \|\mathcal{B}(u^{(n)}, \rho^{(n)} - \rho^{(n-1)})\|_{E_p} + \|\mathcal{B}(u^{(n-1)} - u^{(n)}, \rho^{(n-1)}) \|_{E_p} \nonumber
\\&\le C(\epsilon) \|\rho^{(n-1)} - \rho^{(n-2)} \|_{E_p}
\end{align} where $C(\epsilon)$ is a constant depending on $\epsilon$ obeying $C(\epsilon) < 1$ for a sufficiently small $\epsilon$. Therefore, the sequence $\left\{\rho^{(n)}\right\}_{n=1}^{\infty}$ is Cauchy in $E_p$ and converges to a solution $\rho$ of \eqref{pde}--\eqref{initial}. Uniqueness follows from a similar estimate to \eqref{cauchy}. This finishes the proof of Step 2. Therefore the proof of Theorem~\ref{solinBesov} is complete.

\section{Analyticity of Solutions in Besov Spaces} \la{BesovAnalyt}

In this section, we prove that solutions of \eqref{pde}--\eqref{udef} are analytic in Besov spaces. 

\beg{thm} \la{Analytic} Let $p \in (1, \infty)$. Let $\alpha \in (0, \fr{1}{4}]$. Let $\rho_0 \in \dot{B}_{p,1}^{\fr{2}{p}}(\R^2)$ be sufficiently small. Then the unique solution $\rho \in E_p$ of \eqref{pde}--\eqref{initial}, obtained in Theorem~\ref{solinBesov}, obeys $e^{\alpha t \l_1} \rho \in E_p$ for all $t> 0$, where $\l_1$ is the Fourier multiplier with symbol $|\xi|_1 = |\xi_1| + |\xi_2|$.
\end{thm}

\textbf{Proof:} The main step in the proof is to show that if 
\be 
\rho(t) = e^{-t\l}\rho_0 - \int_{0}^{t} e^{-(t-s)\l}\na \cdot (u\rho)(s) ds,
\ee 
then 
\be 
\|e^{\alpha t \l_1} \rho\|_{E_p} \le C_3\|\rho_0\|_{\dot{B}_{p,1}^{\fr{2}{p}}} + C_4\|e^{\alpha t \l_1} \rho\|_{E_p}^3.
\ee

First, we note that the operator $e^{\alpha t \l_1 - \fr{1}{2}t \l}$ is a Fourier multiplier that is bounded on $L^p$ spaces for $p \in (1, \infty)$. The proof of this latter statement is similar to the proof of Lemma 2 in \cite{BBT}, and is based on the fact that $e^{\alpha t \l_1 - \fr{1}{2}t \l}$ is a Fourier multiplier with symbol 
\be 
e^{\alpha t |\xi|_1 - \fr{1}{2} t|\xi|} \le e^{2\alpha t|\xi| - \fr{1}{2}t |\xi|} = e^{-(\fr{1}{2}-2\alpha)t|\xi|}
\ee which is uniformly bounded by $1$ since $\alpha \in (0, \fr{1}{4}]$. 

Accordingly, for $j \in \mathbb{Z}$, we have 
\be 
\|e^{\alpha t \l_1}e^{-t\l} \Delta_j \rho_0\|_{L^p} = \|e^{\alpha t \l_1 - \fr{1}{2}t\l} e^{-\fr{1}{2}t\l}\Delta_j \rho_0\|_{L^p} \le Ce^{-ct2^j}\|\Delta_j \rho_0\|_{L^p}
\ee 
and so 
\be 
\|e^{\alpha t \l_1}e^{-t\l}\rho\|_{E_p} \le C\|\rho_0\|_{\dot{B}_{p,1}^{\fr{2}{p}}}.
\ee

Now we estimate 
\be 
e^{\alpha t \l_1} \mathcal{B}(u,\rho) = e^{\alpha t \l_1} \int_{0}^{t} e^{-(t-s)\l} \na \cdot (u\rho)(s) ds
\ee in $E_p$.
We start by writing $e^{\alpha t \l_1}\mathcal{B}(u,\rho)$ as 
\be 
e^{\alpha t \l_1} B(u,v) = \int_{0}^{t} e^{\alpha (t-s)\l_1}e^{-\fr{1}{2}(t-s)\l} e^{-\fr{1}{2} (t-s) \l} e^{\alpha s\l_1} \na \cdot (e^{-\alpha s \l_1}\tilde{u} e^{-\alpha s \l_1} \tilde{\rho})(s) ds
\ee 
where 
\be 
\tilde{u}(s) = e^{\alpha s \l_1} u(s)  
\ee and 
\be 
\tilde{\rho}(s) = e^{\alpha s \l_1} \rho(s).
\ee

Using the uniform boundedness of the operator $e^{\alpha (t-s) \l_1 - \fr{1}{2}(t-s) \l}$ on $L^p$ spaces for $p \in (1, \infty)$, Bernstein's inequality, and the bound \eqref{localization}, we get
\be 
\|e^{\alpha t\l_1} \mathcal{B}(u,v)\|_{E_p} \le C\|e^{\alpha s \l_1} \left(e^{-\alpha s \l_1}\tilde{u} e^{-\alpha s \l_1}\tilde{\rho} \right) \|_{\tilde{L}_t^1 \dot{B}_{p,1}^{\fr{2}{p} + 1}}.
\ee
Decomposing $\Delta_j \left(e^{-\alpha s \l_1}\tilde{u} e^{-\alpha s \l_1}\tilde{\rho} \right)$ as 
\beg{align}
\Delta_j \left(e^{-\alpha s \l_1}\tilde{u} e^{-\alpha s \l_1}\tilde{\rho} \right)
&= \sum\limits_{k \ge j-2} \Delta_j \left[\left(e^{-\alpha s \l_1} S_k \tilde{u} \right) \left(e^{-\alpha s \l_1} \Delta_k \tilde{\rho} \right) \right] \nonumber
\\&+ \sum\limits_{k \ge j-2} \Delta_j \left[ \left(e^{-\alpha s \l_1} S_{k+1} \tilde{\rho} \right) \left(e^{-\alpha s \l_1} \Delta_k \tilde{u} \right) \right], 
\end{align}
we have 
\beg{align} 
\left\|e^{\alpha s \l_1} \Delta_j \left(e^{-\alpha s \l_1}\tilde{u} e^{-\alpha s \l_1}\tilde{\rho}  \right) \right\|_{L_t^1 L^p}
&\le C\sum\limits_{k \ge j-2} \left\|e^{\alpha s \l_1} \left(e^{-\alpha s \l_1} S_k \tilde{u} \right) \left(e^{-\alpha s \l_1} \Delta_k \tilde{\rho} \right) \right\|_{L_t^1 L^p} \nonumber
\\&+ C\sum\limits_{k \ge j-2} \left\|e^{\alpha s \l_1} \left(e^{-\alpha s \l_1} S_{k+1} \tilde{\rho} \right) \left(e^{-\alpha s \l_1} \Delta_k \tilde{u} \right) \right\|_{L_t^1 L^p}.
\end{align}

It is shown in \cite{BBT} that the bilinear operator $B_{w}(f,g)$ defined by 
\be \la{BW}
B_{w}(f,g) = e^{w\l_1} (e^{-w \l_1}f e^{-w \l_1}g)
\ee obeys
\be 
\|B_{w}(f,g)\|_{L^p} \le C \|Z_{w}^1 f Z_{w}^2g \|_{L^p}
\ee
where $C > 0$ is a positive constant depending only on $p$, $Z_{w}^1$ and $Z_{w}^2$ are bounded linear operators on $L^p$ for $p \in (1, \infty)$ such that their norms is independent of $w$. For simplicity, we drop the index $w$, and we write $Z^1$ for $Z_w^1$ and $Z^2$ for $Z_w^2$. 

Consequently, 
\beg{align}
\left\|e^{\alpha s \l_1} \Delta_j \left(e^{-\alpha s \l_1}\tilde{u} e^{-\alpha s \l_1}\tilde{\rho}  \right) \right\|_{L_t^1 L^p}
&\le C\sum\limits_{k \ge j-2} \|Z^1 S_k \tilde{u} Z^2 \Delta_k \tilde{\rho}\|_{L_t^1 L^p} \nonumber
\\&+ C\sum\limits_{k \ge j-2} \|Z^1 S_{k+1} \tilde{\rho} Z^2 \Delta_k \tilde{u}\|_{L_t^1 L^p}.  
\end{align}

Now we proceed as in Step 1 of the proof of Theorem \ref{solinBesov}. Indeed, 
\be 
\|Z^1 S_{k+1} \tilde{\rho}\|_{L_t^{\infty} L^{\infty}}
\le \sum\limits_{l \le k} 2^{l\fr{2}{p}} \|Z^1 \Delta_l \tilde{\rho}\|_{L_{t}^{\infty}L^p} 
\le C\sum\limits_{l \le k} 2^{l \fr{2}{p}} \|\Delta_l \tilde{\rho}\|_{L_t^{\infty}L^p}
\le C\|\tilde{\rho}\|_{\tilde{L}_t^{\infty}\dot{B}_{p,1}^{\fr{2}{p}}}.
\ee 
If we show that 
\be \la{Btilde1}
\|\Delta_k \tilde{u}\|_{L_t^1 L^p} \le C \|\tilde{\rho}\|_{L_t^{\infty} \dot{B}_{p,1}^{\fr{2}{p}}} \left(\sum\limits_{m \ge k-2} \|\Delta_m \tilde{\rho}\|_{L_t^1 L^p} \right)
\ee and 
\be \la{Btilde2}
\|Z^1 S_k \tilde{u}\|_{L_t^{\infty}L^{\infty}}  \le C\|\tilde{\rho}\|_{\tilde{L}_t^{\infty}\dot{B}_{p,1}^{\fr{2}{p}}}^2,
\ee
then the rest follows as in Step 1 of Theorem \ref{solinBesov}. 

Hence, we proceed to prove the bounds \eqref{Btilde1} and \eqref{Btilde2}. We note that
\be 
\tilde{u} = e^{\alpha s \l_1} u 
= e^{\alpha s\l_1} \mathbb{P} (\rho R\rho) 
= e^{\alpha s \l_1} \mathbb{P} (e^{-\alpha s \l_1} \tilde{\rho} e^{-\alpha s \l_1} R \tilde{\rho}).
\ee
We decompose $\Delta_l (e^{-\alpha s \l_1} \tilde{\rho} e^{-\alpha s \l_1} R\tilde{\rho})$ as 
\beg{align} 
\Delta_l (e^{-\alpha s \l_1} \tilde{\rho} e^{-\alpha s \l_1} R\tilde{\rho}) 
&= \sum\limits_{m \ge l -2} \Delta_l \left[\left(e^{-\alpha s \l_1} S_{m+1} \tilde{\rho} \right) \left(e^{-\alpha s \l_1}\Delta_m R \tilde{\rho} \right) \right]\nonumber
\\&+ \sum\limits_{m \ge l -2} \Delta_l \left[\left(e^{-\alpha s \l_1} \Delta_m \tilde{\rho} \right) \left(e^{-\alpha s \l_1} S_m R\tilde{\rho} \right)\right].
\end{align}
In view of the boundedness of the operators $Z^1$ and $\mathbb{P}$, we estimate
\beg{align}
&\|Z^1 S_k \tilde{u}\|_{L_t^{\infty}L^{\infty}}
\le C \sum\limits_{l \le k-1} 2^{l \fr{2}{p}} \|\Delta_l \tilde{u}\|_{L_t^{\infty}L^p} \nonumber
\le C \sum\limits_{l \le k-1} 2^{l \fr{2}{p}} \|\Delta_l \left[e^{\alpha s \l_1} (e^{-\alpha s \l_1} \tilde{\rho} e^{-\alpha s \l_1} R\tilde{\rho}) \right]\|_{L_t^{\infty} L^p} \nonumber
\\&\le C \sum\limits_{l \le k-1} \sum\limits_{m \ge l-2} 2^{l\fr{2}{p}} \|B_{\alpha s}(S_{m+1} \tilde{\rho}, \Delta_m R 
\tilde{\rho})\|_{L_t^{\infty}L^p} \nonumber
+ C\sum\limits_{l \le k-1} \sum\limits_{m \ge l-2} 2^{l\fr{2}{p}} \|B_{\alpha s}(S_m R\tilde{\rho}, \Delta_m 
\tilde{\rho})\|_{L_t^{\infty}L^p}
\end{align} where $B_{w}(f,g)$ is defined in \eqref{BW}. This implies that 
\beg{align}
\|Z^1 S_k \tilde{u}\|_{L_t^{\infty}L^{\infty}}
&\le  C \sum\limits_{l \le k-1} \sum\limits_{m \ge l-2} 2^{l\fr{2}{p}} \|(Z^1 S_{m+1} \tilde{\rho}) (Z^2 \Delta_m R 
\tilde{\rho})\|_{L_t^{\infty}L^p} \nonumber
\\&+ C\sum\limits_{l \le k-1} \sum\limits_{m \ge l-2} 2^{l\fr{2}{p}} \|(Z^1 S_m R\tilde{\rho}) (Z^2 \Delta_m 
\tilde{\rho})\|_{L_t^{\infty}L^p}.
\end{align} Now we proceed as in \eqref{N1} and \eqref{N2} and we obtain \eqref{Btilde2}.
Finally, we estimate
\beg{align}
&\|\Delta_k \tilde{u}\|_{L_t^1 L_p}
\le C\|\Delta_k e^{\alpha s \l_1} \left[(e^{-\alpha s \l_1}\tilde{\rho}) (e^{-\alpha s \l_1} R \tilde{\rho})  \right]\|_{L_t^1 L^p} \nonumber
\\&\le C \sum\limits_{m \ge k-2} \|e^{\alpha s \l_1} \left[(e^{-\alpha s \l_1} S_{m+1} \tilde{\rho}) (e^{-\alpha s \l_1} \Delta_m R\tilde{\rho}) \right] \|_{L_t^1 L^p} \nonumber
\\&+ C \sum\limits_{m \ge k-2} \|e^{\alpha s \l_1} \left[(e^{-\alpha s \l_1} S_m R\tilde{\rho}) (e^{-\alpha s \l_1} \Delta_m \tilde{\rho}) \right] \|_{L_t^1 L^p} \nonumber
\\&= C \sum\limits_{m \ge k-2} \|B_{\alpha s} (S_{m+1} \tilde{\rho}, \Delta_m R\tilde{\rho}) \|_{L_t^1 L^p} \nonumber
+ C \sum\limits_{m \ge k-2} \|B_{\alpha s} (S_m R\tilde{\rho}, \Delta_m \tilde{\rho}) \|_{L_t^1 L^p} \nonumber
\\&\le C \sum\limits_{m \ge k-2} \|Z^1S_{m+1} \tilde{\rho}\|_{L_t^{\infty}L^{\infty}} \|Z^2\Delta_m R \tilde{\rho}\|_{L_t^1 L^p} \nonumber
+C \sum\limits_{m \ge k-2} \|Z^1 S_{m} R \tilde{\rho}\|_{L_t^{\infty}L^{\infty}} \|Z^2\Delta_m \tilde{\rho}\|_{L_t^1 L^p}
\\&\le C\|\tilde{\rho}\|_{\tilde{L}_t^{\infty} \dot{B}_{p,1}^{\fr{2}{p}}} \left(\sum\limits_{m \ge k-2} \|\Delta_m \tilde{\rho}\|_{L_t^1 L^p} \right)
\end{align}
which proves \eqref{Btilde1}. This ends the proof of Theorem \ref{Analytic}.

\section{Regularity of Solutions for Small Initial Data} \la{reg}

In this section, we study the regularity of solutions of \eqref{pde}--\eqref{initial} for small initial data. 

We use the following  lemma:

\beg{lem} \la{lem1} Let $j \in \mathbb{Z}$, $t > 0$, $\alpha \in [0,1), c>0.$ Then there is a positive constant $C > 0$ that depends on $\alpha$ but does not dependent of $j$ and $t$ such that the estimate 
\be \la{ine}
\int_{0}^{t} 2^j e^{-c(t-s)2^j} s^{-\alpha} ds \le Ct^{-\alpha}
\ee holds. 
\end{lem}

\textbf{Proof:} We split the given integral into the sum 
\be 
\int_{0}^{t} 2^j e^{-c(t-s)2^j} s^{-\alpha} ds 
= \int_{0}^{\fr{t}{2}} 2^j e^{-c(t-s)2^j} s^{-\alpha} ds
+ \int_{\fr{t}{2}}^{t} 2^j e^{-c(t-s)2^j} s^{-\alpha} ds.
\ee Using the fact that $2^{j}e^{-c(t-s)2^j} \le C(t-s)^{-1}$ for all $s \in [0, \fr{t}{2}]$, we estimate
\be \la{ine1}
\int_{0}^{\fr{t}{2}} 2^j e^{-c(t-s)2^j} s^{-\alpha} ds 
\le C\int_{0}^{\fr{t}{2}} (t-s)^{-1} s^{-\alpha} ds
\le Ct^{-1} \int_{0}^{\fr{t}{2}} s^{-\alpha} ds
\le C_{\alpha} t^{-\alpha}.
\ee Using the fact that $s^{-\alpha} \le 2^{\alpha} t^{-\alpha}$ for all $s \in [\fr{t}{2}, t]$, we estimate
\be \la{ine2}
\int_{\fr{t}{2}}^{t} 2^j e^{-c(t-s)2^j} s^{-\alpha} ds
\le C_{\alpha} t^{-\alpha} \int_{\fr{t}{2}}^{t} 2^j e^{-c(t-s)2^j} ds
= C_{\alpha}t^{-\alpha} \left[1 - e^{-2^{j-1}ct} \right]
\le C_{\alpha} t^{-\alpha}.
\ee
Adding \eqref{ine1} and \eqref{ine2}, we obtain \eqref{ine}.

\beg{thm} \la{hreg} Let $\alpha \in [0,1), \beta > 0$. Let $\rho_0 \in \dot{B}_{2,1}^{1}(\R^2) \cap \dot{B}_{\infty, \infty}^{\beta - \alpha}(\R^2)$ be  sufficiently small. Then there is a positive constant $C> 0$ depending on the initial data such that  the unique solution $\rho$ of \eqref{pde}--\eqref{initial} satisfies 
\be \la{highreg}
\sup\limits_{t > 0} t^{\alpha} \|\rho(t)\|_{\dot{B}_{\infty, \infty}^{\beta}}
\le C.
\ee
\end{thm}

\textbf{Proof:} We consider the approximating initial value problem \eqref{approxeq} whose solution is given by 
\be 
\rho^{(n)}(t) = e^{-t\l}\rho_0 - \mathcal{B}(u^{(n-1)}, \rho^{(n-1)})
\ee 
where $\mathcal{B}$ is the bilinear form defined by 
\be 
\mathcal{B}(v, \theta) = \int_{0}^{t} e^{-(t-s)\l} \na \cdot (v \theta)(s) ds. 
\ee

First, we estimate
\be 
t^{\alpha}2^{j\beta} \|e^{-t\l} \Delta_j \rho_0\|_{L^{\infty}}
\le Ct^{\alpha}e^{-ct2^j} 2^{j\beta}\|\Delta_j \rho_0\|_{L^{\infty}}
\le C 2^{-j \alpha} 2^{j \beta}\|\Delta_j \rho_0\|_{L^{\infty}}
\le C \|\rho_0\|_{\dot{B}_{\infty, \infty}^{\beta - \alpha}} 
\ee in view of \eqref{localization} and the bound $x^{\alpha}e^{-x} \le C$ that holds for all $x \ge 0$. 

We show that \be \la{RTP}
\sup\limits_{t > 0} \left\{t^{\alpha}\|\mathcal{B}(u^{(n-1)}, \rho^{(n-1)})(t)\|_{\dot{B}_{\infty, \infty}^{\beta}} \right\}
\le C \|\rho^{(n-1)}\|_{\tilde{L}_t^{\infty} \dot{B}_{p,1}^{\fr{2}{p}}}^2 \sup\limits_{t > 0} \left\{t^{\alpha} \|\rho^{(n-1)}(t)\|_{\dot{B}_{\infty, \infty}^{\beta}}\right\}.
\ee 
We start by applying $\Delta_j$ to $\mathcal{B}(u^{(n-1)}, \rho^{n-1)})$, we use the paraproduct decomposition 
\be 
\Delta_j (u^{(n-1)} \rho^{(n-1)}) 
= \sum\limits_{k \ge j-2} \Delta_j (S_k u^{(n-1)} \Delta_k \rho^{(n-1)})
+ \sum\limits_{k \ge j-2} \Delta_j (S_{k+1} \rho^{(n-1)} \Delta_k u^{(n-1)}),
\ee and we obtain 
\be 
\Delta_j \mathcal{B}(u^{(n-1)}, \rho^{(n-1)}) = \mathcal{B}_{1,j}(u^{(n-1)}, \rho^{(n-1)}) + \mathcal{B}_{2,j}(u^{(n-1)}, \rho^{(n-1)}) 
\ee where 
\be 
\mathcal{B}_{1,j}(u^{(n-1)}, \rho^{(n-1)}) = \int_{0}^{t} e^{-(t-s)\l}\na \cdot \left[\sum\limits_{k \ge j-2} \Delta_j (S_k u^{(n-1)} \Delta_k \rho^{(n-1)}) \right](s) ds,
\ee 
and
\be 
\mathcal{B}_{2,j}(u^{(n-1)}, \rho^{(n-1)})= \int_{0}^{t} e^{-(t-s)\l} \na \cdot \left[\sum\limits_{k \ge j-2} \Delta_j (S_{k+1} \rho^{(n-1)} \Delta_k u^{(n-1)}) \right](s) ds.
\ee
In view of Bernstein's inequality \eqref{BER1}, the bounds \eqref{localization} and  \eqref{nonlinear1}, and Lemma \ref{lem1}, we estimate
\beg{align}
&2^{j\beta} \|\mathcal{B}_{1,j}(u^{(n-1)}, \rho^{(n-1)})\|_{L^{\infty}}
\le C 2^{j\beta} 2^{j} \int_{0}^{t} e^{-c(t-s)2^j} \left[\sum\limits_{k \ge j-2} \|S_{k} u^{(n-1)}\|_{L^{\infty}} \|\Delta_k \rho^{(n-1)}\|_{L^{\infty}} \right] ds \nonumber
\\&\le C\|\rho^{(n-1)}\|_{\tilde{L}_{t}^{\infty} \dot{B}_{2,1}^{1}}^2 \left\{\int_{0}^{t} \left(2^j e^{-c(t-s)2^j}s^{-\alpha} \right) \left(\sum\limits_{k \ge j-2} 2^{- (k-j)\beta} s^{\alpha} 2^{k\beta} \|\Delta_k \rho^{(n-1)}\|_{L^{\infty}} \right)    ds \right\} \nonumber
\\&\le C\|\rho^{(n-1)}\|_{\tilde{L}_{t}^{\infty} \dot{B}_{2,1}^{1}}^2 \sup\limits_{t>0} \left\{t^{\alpha} \|\rho^{(n-1)}\|_{\dot{B}_{\infty, \infty}^{\beta}} \right\} \int_{0}^{t} 2^j e^{-c(t-s)2^j}s^{-\alpha} ds  \nonumber
\\&\le Ct^{-\alpha} \|\rho^{(n-1)}\|_{\tilde{L}_{t}^{\infty} \dot{B}_{2,1}^{1}}^2 \sup\limits_{t>0} \left\{t^{\alpha} \|\rho^{(n-1)}\|_{\dot{B}_{\infty, \infty}^{\beta}} \right\},
\end{align} hence 
\beg{align} \la{B1BOUND}
t^{\alpha} 2^{j\beta} \|\mathcal{B}_{1,j}(u^{(n-1)}, \rho^{(n-1)})\|_{L^{\infty}} \le  C \|\rho^{(n-1)}\|_{\tilde{L}_{t}^{\infty} \dot{B}_{2,1}^{1}}^2 \sup\limits_{t>0} \left\{t^{\alpha} \|\rho^{(n-1)}\|_{\dot{B}_{\infty, \infty}^{\beta}} \right\}.
\end{align}

Now, we estimate $2^{j\beta} \|\mathcal{B}_{2,j}(u^{(n-1)}, \rho^{(n-1)})\|_{L^{\infty}}$. We note first that 
\beg{align}
&\sum\limits_{k \ge j-2} s^{\alpha} 2^{j \beta} \|\Delta_k u^{(n-1)}\|_{L^{\infty}} \nonumber
\le C \sum\limits_{k \ge j-2} s^{\alpha}  2^{j \beta} \|\Delta_k (\rho^{(n-1)}R\rho^{(n-1)})\|_{L^{\infty}} \nonumber
\\&\le C \sum\limits_{k \ge j-2} s^{\alpha}  2^{j \beta} \left(\sum\limits_{l\ge k-2} \|S_{l+1} \rho^{(n-1)} \|_{L^{\infty}} \|\Delta_l R\rho^{(n-1)}\|_{L^{\infty}} \right) \nonumber
\\&+  C \sum\limits_{k \ge j-2}  s^{\alpha}  2^{j \beta} \left(\sum\limits_{l\ge k-2} \|S_l R\rho^{(n-1)}\|_{L^{\infty}} \|\Delta_l \rho^{(n-1)}\|_{L^{\infty}} \right)  \nonumber
\\&\le C \|\rho^{(n-1)}\|_{\tilde{L}_t^{\infty} \dot{B}_{2,1}^{1}} \sum\limits_{k \ge j-2} \sum\limits_{l \ge k-2} 2^{-(l-j)\beta} s^{\alpha}  2^{l \beta} \|\Delta_l \rho^{(n-1)}\|_{L^{\infty}} \nonumber
\\&= C \|\rho^{(n-1)}\|_{\tilde{L}_t^{\infty} \dot{B}_{2,1}^{1}} \sum\limits_{l \ge j-4} \sum\limits_{j-2 \le k \le l+2} 2^{-(l-j)\beta} s^{\alpha}  2^{l \beta} \|\Delta_l \rho^{(n-1)}\|_{L^{\infty}} \nonumber
\\&= C \|\rho^{(n-1)}\|_{\tilde{L}_t^{\infty} \dot{B}_{2,1}^{1}} \sum\limits_{l \ge j-4} (l-j+5) 2^{-(l-j)\beta} s^{\alpha}  2^{l \beta} \|\Delta_l \rho^{(n-1)}\|_{L^{\infty}} \nonumber
\\&\le C \|\rho^{(n-1)}\|_{\tilde{L}_t^{\infty} \dot{B}_{2,1}^{1}} \sup\limits_{t>0} \left\{t^{\alpha} \|\rho^{(n-1)}\|_{\dot{B}_{\infty, \infty}^{\beta}} \right\}.
\end{align}
Here, we have used the boundedness of the Leray projector in Besov spaces, the paraproduct decomposition \eqref{para}, the bound \eqref{nonlinear1}, Fubini's theorem, and Young's convolution inequality. This implies 
\beg{align}
2^{j\beta} \|\mathcal{B}_{2,j}(u^{(n-1)}, \rho^{(n-1)})\|_{L^{\infty}}
&\le C 2^{j\beta} 2^{j} \int_{0}^{t} e^{-c(t-s)2^j} \left[\sum\limits_{k \ge j-2} \|S_{k+1} \rho^{(n-1)}\|_{L^{\infty}} \|\Delta_k u^{(n-1)}\|_{L^{\infty}}\right] ds \nonumber
\\&\le C\|\rho^{(n-1)}\|_{\tilde{L}_{t}^{\infty} \dot{B}_{2,1}^{1}}^2 \sup\limits_{t>0} \left\{t^{\alpha} \|\rho^{(n-1)}\|_{\dot{B}_{\infty, \infty}^{\beta}} \right\} \int_{0}^{t} 2^j e^{-c(t-s)2^j}s^{-\alpha} ds  \nonumber
\\&\le Ct^{-\alpha} \|\rho^{(n-1)}\|_{\tilde{L}_{t}^{\infty} \dot{B}_{2,1}^{1}}^2 \sup\limits_{t>0} \left\{t^{\alpha} \|\rho^{(n-1)}\|_{\dot{B}_{\infty, \infty}^{\beta}} \right\},
\end{align} hence 
\beg{align} \la{B2BOUND}
t^{\alpha} 2^{j\beta} \|\mathcal{B}_{2,j}(u^{(n-1)}, \rho^{(n-1)})\|_{L^{\infty}} \le  C \|\rho^{(n-1)}\|_{\tilde{L}_{t}^{\infty} \dot{B}_{2,1}^{1}}^2 \sup\limits_{t>0} \left\{t^{\alpha} \|\rho^{(n-1)}\|_{\dot{B}_{\infty, \infty}^{\beta}} \right\}.
\end{align}

Putting \eqref{B1BOUND} and \eqref{B2BOUND} together, we obtain \eqref{RTP}.
Therefore, 
\be 
\sup\limits_{t > 0} t^{\alpha} \|\rho^{(n)}(t)\|_{\dot{B}_{\infty, \infty}^{\beta}} \le C_3\|\rho_0\|_{\dot{B}_{\infty, \infty}^{\beta - \alpha}} + C_4 \|\rho^{(n-1)}\|_{\tilde{L}_{t}^{\infty} \dot{B}_{2,1}^{1}}^2 \sup\limits_{t>0} \left\{t^{\alpha} \|\rho^{(n-1)}(t)\|_{\dot{B}_{\infty, \infty}^{\beta}} \right\}.
\ee
Now, we use the smallness of the initial data and we proceed as in Theorem \ref{solinBesov}. We omit further details. 

We recall the following relationship between the inhomogeneous Besov space $B_{\infty, \infty}^s$ and  H\"older spaces: 

\beg{rem} \la{nonhomreg}
For $s \in \R^{+} \setminus \N$, the inhomogeneous Besov space $B_{\infty, \infty}^s(\R^2)$ coincides with the H\"older space $C^{[s], s-[s]}(\R^2)$ of bounded functions $f$ whose derivatives of order $|\alpha| \le s$ are bounded and satisfy 
\be 
|\pa^{\alpha} f(x) - \pa^{\alpha}f(y)| \le C|x-y|^{s-[s]}
\ee for $|x-y| \le 1$ (see \cite{BCD}).
\end{rem}

As a consequence, we obtain the following regularity result: 

\beg{cor} \la{last} Let $s \ge 1$ be an integer. Let $\rho_0 \in L^{\infty}(\R^2) \cap \dot{B}_{2,1}^{1}(\R^2) \cap \dot{B}_{\infty, \infty}^{s + \fr{1}{2}}(\R^2)$ be sufficiently small. Let $\rho$ be the unique solution of \eqref{pde}--\eqref{initial}. Then $D^{\gamma} \rho \in L^{\infty}(\R^2)$ for $|\gamma| \le s$, and its $L^{\infty}$ norm is uniformly bounded in time. Moreover, for $|\gamma| \le s$, $D^{\gamma}\rho$ is H\"older continuous with a uniform in time H\"older bound.  
\end{cor}

\textbf{Proof:} In view of \eqref{hominhom}, the bound \eqref{decayinf}, and Theorem \ref{hreg} applied with $\alpha = 0$ and $\beta = s + \fr{1}{2}$, we have 
\be 
\|\rho(t)\|_{B_{\infty, \infty}^{s+ \fr{1}{2}}} 
\le C \left\{\|\rho(t)\|_{\dot{B}_{\infty, \infty}^{s+ \fr{1}{2}}} + \|\rho(t)\|_{L^{\infty}} \right\} 
\le C (1 + \|\rho_0\|_{L^{\infty}})
\ee where $C$ is a constant depending only on the initial data. Remark \ref{nonhomreg} completes the proof of Corollary  \ref{last}. 

We consider the long time behavior of derivatives of solutions of \eqref{pde}--\eqref{initial} for sufficiently small initial data in Besov spaces:

\beg{cor} Let $s \ge 1$ be an integer. Let $\delta \in (0,1)$. Let $\rho_0 \in L^{\infty}(\R^2) \cap \dot{B}_{2,1}^{1}(\R^2) \cap \dot{B}_{\infty, \infty}^{s + \delta}(\R^2)$ be sufficiently small. Let $\rho$ be the unique solution of \eqref{pde}--\eqref{initial}. Then 
\be \la{longtime}
\lim\limits_{t \to \infty} \left\{\|D^{\gamma} \rho(t)\|_{L^{\infty}} + [D^{\gamma}  \rho(t)]_{\delta} \right\} = 0 
\ee for all $|\gamma| \le s$, where 
\be 
[D^{\gamma} \rho(t)]_{\delta} = \sup\limits_{0<|x-y| \le 1} \fr{|D^{\gamma}\rho(y) - D^{\gamma} \rho(x)|}{|x-y|^{\delta}}. 
\ee
\end{cor}

\textbf{Proof:} We show that $\rho(\cdot,t) \in H^{2}(\R^2)$ in order to apply Remark~\ref{ltime}. Indeed,
\beg{align} 
&\int_{0}^{\infty} \|\rho(t)\|_{\dot{B}_{2,2}^2} dt 
\le C\int_{0}^{\infty} \|\rho(t)\|_{\dot{B}_{2,1}^2}dt 
= C\int_{0}^{\infty} \sum\limits_{j \in \ZZ} 2^{2j} \|\Delta_j \rho(t)\|_{L^2} dt \nonumber
\\&= C \sum\limits_{j \in \ZZ} 2^{2j} \|\Delta_j \rho(t)\|_{L_t^1 L^2} 
= C \|\rho\|_{\tilde{L}_t^1 \dot{B}_{2,1}^2} < \infty
\end{align} in view of the continuous Besov embedding \eqref{BESOVEMBED} and the monotone convergence theorem. 
But $B_{2,2}^{2}$ coincides with the Sobolev space $H^2$. Thus, 
\be 
 \|\rho (t) \|_{H^2} < \infty 
\ee for a.e. $t \in (0, \infty)$, and so $\rho(\cdot,t) \in H^2$ for a.e. $t \in [0,\infty)$.

In view of Remark \ref{nonhomreg}, Remark \ref{ltime}, and Theorem \ref{hreg} applied with $\alpha = \fr{1}{2}$ and $\beta = s + \delta$, we obtain    
\beg{align}
&\|D^{\gamma} \rho(t)\|_{L^{\infty}} + [D^{\gamma}  \rho(t)]_{\delta} 
\leq C \|\rho(t)\|_{B_{\infty, \infty}^{s+ \delta}} \nonumber
\\&\le C \left\{\|\rho(t)\|_{\dot{B}_{\infty, \infty}^{s+ \delta}} + \|\rho(t)\|_{L^{\infty}} \right\} 
\le C \left(\fr{1}{\sqrt t} + \fr{\|\rho_0\|_{L^{\infty}}}{1 + Ct\|\rho_0\|_{L^{\infty}}}\right).  
\end{align}
Letting $t \rightarrow \infty$, we obtain \eqref{longtime}.

\section{Regularity of Solutions for Arbitrary Initial Data} \la{regarb}

In this section, we prove that any solution of \eqref{pde}--\eqref{initial} is smooth for arbitrary initial data, provided that it satisfies a certain regularity condition. 

\beg{thm} \la{arbitthm} Let $\rho$ be a weak solution of \eqref{pde}--\eqref{initial} on $[0, \infty)$.
Let $0 < t_0 < t < \infty$. If
\be \la{arbitinitialcond}
\rho \in L^{\infty} ([t_0, t]; C^{\delta}(\R^2)),
\ee 
for some $\delta \in (0,1)$, then 
\be 
\rho \in C^{\infty}((t_0, t] \times \R^2).
\ee 
\end{thm} 

\textbf{Proof:} 
We sketch the main ideas. 
Let us note first that 
\be  \la{arbit2}
u \in L^{\infty} ([t_0, t]; C^{\delta}(\R^2)).
\ee 
where 
\be 
u = - \PP(\rho R \rho).
\ee
Indeed, for any $s \in [t_0, t]$, we have 
\beg{align} 
\|u(s)\|_{C^{\delta}} 
&\le C\|\rho (s) R \rho (s)\|_{C^{\delta}} \nonumber
\\&\le C\|\rho (s)\|_{L^{\infty}} \|R\rho (s)\|_{L^{\infty}} 
+ C \|\rho (s)\|_{L^{\infty}} \|R\rho (s)\|_{C^{\delta}}
+ C\|R\rho (s) \|_{L^{\infty}} \|\rho (s)\|_{C^{\delta}} \nonumber
\\&\le C \|\rho (s)\|_{C^{\delta}}^2
\end{align} in view of the boundedness of the Leray projector and Riesz transforms on the H\"older space $C^{\delta}$. Consequently, the H\"older regularity of $\rho$ imposed in \eqref{arbitinitialcond} gives \eqref{arbit2}.

Next, we show that 
\be \la{arbit3}
\rho \in L^{\infty} ([t_0, t]; \dot{B}_{p, \infty}^{\delta_1}(\R^2) \cap C^{\delta_1}(\R^2)) 
\ee 
and 
\be \la{arbit4}
u \in L^{\infty} ([t_0, t]; \dot{B}_{p, \infty}^{\delta_1}(\R^2) \cap C^{\delta_1}(\R^2))
\ee for any $p \ge 2$ and $\delta_1 = \delta \left(1 - \fr{2}{p} \right)$.
Indeed, for any $s \in [t_0, t]$, we have
\beg{align}
\|\rho (s)\|_{\dot{B}_{p,\infty}^{\delta_1}} 
&= \sup\limits_{j \in \ZZ} \left( 2^{\delta_1 j} \|\Delta_j \rho (s)\|_{L^p} \right)
\le \sup\limits_{j \in \ZZ} \left( 2^{\delta_1 j} \|\Delta_j \rho (s)\|_{L^{\infty}}^{1- \fr{2}{p}} \|\Delta_j \rho (s)\|_{L^2}^{\fr{2}{p}} \right) \nonumber
\\&\le C \left(\|\rho (s)\|_{\dot{B}_{\infty, \infty}^{\delta}}\right)^{1-\fr{2}{p}} \|\rho(s)\|_{L^2}^{\fr{2}{p}}
\le C\left(\|\rho (s)\|_{C^{\delta}}\right)^{1-\fr{2}{p}} \|\rho(s)\|_{L^2}^{\fr{2}{p}}
\end{align} and similarly
\beg{align}
\|u (s)\|_{\dot{B}_{p,\infty}^{\delta_1}} 
\le C \left(\|u (s)\|_{\dot{B}_{\infty, \infty}^{\delta}}\right)^{1-\fr{2}{p}} \|u(s)\|_{L^2}^{\fr{2}{p}}
\le C\left(\|u (s)\|_{C^{\delta}}\right)^{1-\fr{2}{p}} \|\rho(s)\|_{L^4}^{\fr{4}{p}}.
\end{align} 
The last inequality holds in view of the boundedness of the Leray projector on $L^2$ followed by an application of H\"older's inequality with exponents $4,4$. 
The interpolation inequality 
\be 
\|\rho(s)\|_{L^4} \le \|\rho(s)\|_{L^{\infty}}^{1/2} \|\rho(s)\|_{L^2}^{1/2}
\ee
together with \eqref{arbit2} and \eqref{arbitinitialcond} gives \eqref{arbit3} and \eqref{arbit4}.

Now, we proceed as in \cite{CW}. We apply $\Delta_j$ to \eqref{pde}, we multiply the resulting equation by $p|\Delta_j \rho|^{p-2} \Delta_j \rho$, we integrate first in the space variable $x \in \R^2$ and then in time from $t_0$ to $t$. We obtain the bound
\beg{align} 
&\|\Delta_j \rho (t)\|_{L^p}
\le Ce^{-c2^j(t-t_0)}\|\Delta_j \rho (t_0)\|_{L^p} \nonumber
\\&+ C\int_{t_0}^{t} e^{-c2^j (t-s)}2^{(1-2\delta_1)j} \left(\|\rho (s)\|_{C^{\delta_1}} \|u(s)\|_{\dot{B}_{p, \infty}^{\delta_1}} + \|u(s)\|_{C^{\delta_1}} \|\rho(s)\|_{\dot{B}_{p, \infty}^{\delta_1}} \right) ds
\end{align} (see \cite{CW} for details). We multiply by $2^{2\delta_1 j}$ and we take the $\ell^{\infty}$ norm in $j$. This yields the bound
\beg{align}
&\|\rho(t)\|_{\dot{B}_{p, \infty}^{2\delta_1}} 
\le C \sup\limits_{j \in \ZZ} \left\{2^{\delta_1 j} e^{-c2^j (t-t_0)} \right\} \|\rho(t_0)\|_{\dot{B}_{p, \infty}^{\delta_1}} \nonumber
\\&+ C \sup\limits_{j \in \ZZ} \left\{1 - e^{-c2^j (t-t_0)} \right\} \sup_{s \in [t_0, t]}  \left\{ \|\rho (s)\|_{C^{\delta_1}} \|u(s)\|_{\dot{B}_{p, \infty}^{\delta_1}} + \|u(s)\|_{C^{\delta_1}} \|\rho(s)\|_{\dot{B}_{p, \infty}^{\delta_1}}\right\}.
\end{align} 
Therefore 
\be 
\rho(\cdot,t) \in \dot{B}_{p, \infty}^{2\delta_1}(\R^2).
\ee 
for any $p \geq 2$. 
In view of the continuous Besov embedding \eqref{BESOVEMBED}, we have the continuous inclusion
\be 
\dot{B}_{p, \infty}^{2\delta_1}(\R^2)\hookrightarrow \dot{B}_{\infty, \infty}^{2\delta_1 - \fr{2}{p}} (\R^2) 
\ee for any $p \geq 2$. We choose $p > \fr{2+2\delta}{\delta}$ so that $2\delta_1 - \fr{2}{p} > \delta_1$, hence 
\be 
\rho(\cdot,t) \in \dot{B}_{p, \infty}^{\delta_2}(\R^2) \cap C^{\delta_2}(\R^2)
\ee 
where $\delta_2 > \delta_1$. We iterate the above process finitely many times and we obtain 
\be 
\rho(\cdot, t) \in C^{\gamma}(\R^2)
\ee 
for some $\gamma> 1$. This ends the proof of Theorem \ref{arbitthm}.

\section{Periodic case} \la{per}

In this section, we consider the initial value problem \eqref{pde}--\eqref{initial} posed on the torus $\mathbb{T}^2$ with periodic boundary conditions. We assume the initial data $\rho_0$ have zero mean. 
We prove existence and regularity of solutions.

\beg{thm} \la{solinBesovper} Let $1 \le p < \infty$. Let $\rho_0 \in \dot{B}_{p,1}^{\fr{2}{p}}(\mathbb{T}^2)$ be sufficiently small. We consider the functional space $E_p$ defined by 
\be 
E_p (\mathbb{T}^2) = \left\{f (t) \in \mathcal{D}_0'(\mathbb{T}^2) : \|f\|_{E_p(\mathbb{T}^2)} = \|f\|_{\tilde{L}_t^{\infty} \dot{B}_{p,1}^{\fr{2}{p}}(\mathbb{T}^2)} + \|f\|_{\tilde{L}_t^{1} \dot{B}_{p,1}^{\fr{2}{p} +1 }(\mathbb{T}^2)} < \infty \right\}.
\ee 
Then \eqref{pde}--\eqref{initial} has a unique global in time solution $\rho \in E_p(\mathbb{T}^2)$.
\end{thm}

The proof of Theorem \ref{solinBesovper} follows from the proof of Theorem \ref{solinBesov}.

In view of the Besov embedding and Theorem \ref{solinBesovper}, we conclude that if $\rho_0 \in \dot{B}_{2,1}^{1}(\mathbb{T}^2)$ is sufficiently small, then there is a constant $C > 0$ depending only on the initial data such that the unique solution $\rho$ of \eqref{pde}--\eqref{initial} obeys
\be \la{uniformboundsol}
\sup\limits_{t > 0} \|\na \rho (t)\|_{L^2(\mathbb{T}^2)} + \int_{0}^{\infty} \|\Delta \rho (t)\|_{L^2(\mathbb{T}^2)} dt \le C. 
\ee  

Using this latter estimate, we end this section by showing that the $L^2(\mathbb{T}^2)$ norm of $\l^{\fr{1}{2}} \rho$ converges exponentially in time to zero. 

We use the following uniform Gronwall lemma \cite{AI}:

\begin{lem} \la{exp} Let $y(t) \geq 0$ obey a differential inequality $$\frac{d}{dt}y + c_1y \leq F_1 + F(t)$$ with initial datum $y(0) = y_0$ with $F_1$ a nonnegative constant, and $F(t) \geq 0$ obeying $$\int\limits_{t}^{t+1} F(s) ds \leq g_0e^{-c_2t} + F_2$$ where $c_1, c_2, g_0$ are positive constants and $F_2$ is a nonnegative constant. Then $$y(t) \leq y_0e^{-c_1t} + g_0e^{c_1+c}(t+1)e^{-ct} + \frac{1}{c_1}F_1 + \frac{e^{c_1}}{1-e^{-c_1}} F_2$$ holds with $c = \min \left\{c_1, c_2 \right\}$.
\end{lem}

\beg{cor}
Let $\rho_0 \in \dot{B}_{2,1}^1 (\mathbb{T}^2)$ be sufficiently small. Then there is a constant $C>0$ depending only on the initial data such that the unique solution $\rho$ of \eqref{pde}--\eqref{initial} obeys 
\be \la{decayhhalf}
\|\l^{\fr{1}{2}} \rho(t)\|_{L^2(\mathbb{T}^2)}^2 \le Ce^{-t} 
\ee for all $t \ge 0$. 
\end{cor}

\textbf{Proof:} We take the inner product in $L^2(\mathbb{T}^2)$ of \eqref{pde} with $\l \rho$ to obtain
\be 
\fr{1}{2} \fr{d}{dt} \|\l^{\fr{1}{2}} \rho(t)\|_{L^2(\mathbb{T}^2)}^2 + \|\l \rho(t)\|_{L^2(\mathbb{T}^2)}^2 = - \int_{\mathbb{T}^2} (u \cdot \na \rho) \l \rho dx.
\ee 
We estimate the nonlinear term 
\beg{align} 
\left|\int_{\mathbb{T}^2} (u \cdot \na \rho) \l \rho dx \right| 
&\le C\|\rho\|_{L^{\infty}(\mathbb{T}^2)}\|\rho\|_{L^4(\mathbb{T}^2)} \|\na \rho\|_{L^4(\mathbb{T}^2)}\|\na \rho\|_{L^2(\mathbb{T}^2)} \nonumber
\\&\le C\|\rho\|_{L^4(\mathbb{T}^2)} \|\na \rho\|_{L^4(\mathbb{T}^2)}^2 \|\na \rho\|_{L^2(\mathbb{T}^2)} \nonumber
\\&\le C\|\rho\|_{L^2(\mathbb{T}^2)}^{\fr{1}{2}} \|\na \rho\|_{L^2(\mathbb{T}^2)}^{\fr{5}{2}} \|\Delta \rho\|_{L^2(\mathbb{T}^2)}
\end{align}
in view of the boundedness of the Leray projector and Riesz transforms on $L^4(\mathbb{T}^2)$, the continuous embedding $W^{1,4}(\mathbb{T}^2)\hookrightarrow L^{\infty}(\mathbb{T}^2)$, and the Ladyzhenskaya interpolation inequality. 

Since $H^{1}(\mathbb{T}^2)$ is continuously embedded in $H^{\fr{1}{2}}(\mathbb{T}^2)$, we have 
\be 
\|\l^{\fr{1}{2}} \rho\|_{L^2(\mathbb{T}^2)} \le C\|\l \rho\|_{L^2(\mathbb{T}^2)},
\ee yielding the differential inequality 
\be \la{exphhalf}
\fr{d}{dt} \|\l^{\fr{1}{2}} \rho\|_{L^2(\mathbb{T}^2)} + C_1 \|\l^{\fr{1}{2}}\rho \|_{L^2(\mathbb{T}^2)} 
\le C_2 \|\rho\|_{L^2(\mathbb{T}^2)}^{\fr{1}{2}} \|\na \rho\|_{L^2(\mathbb{T}^2)}^{\fr{5}{2}} \|\Delta \rho\|_{L^2(\mathbb{T}^2)}.
\ee

We note that \be \la{expl2}
\|\rho(t)\|_{L^2(\mathbb{T}^2)} \le C\|\rho_0\|_{L^2(\mathbb{T}^2)}e^{-ct}
\ee for all $t \ge 0$. Indeed, we multiply \eqref{pde} by $\rho$ and we integrate in the space variable. Then we use the cancellation of the nonlinear term and the continuous embedding of $H^{\fr{1}{2}}(\mathbb{T}^2)$ in $L^2(\mathbb{T}^2)$ to obtain
\be 
\fr{d}{dt} \|\rho(t)\|_{L^2(\mathbb{T}^2)} + C \|\rho(t)\|_{L^2(\mathbb{T}^2)} \le 0
\ee which gives \eqref{expl2}.

Now we go back to the differential inequality \eqref{exphhalf}. Using the bounds \eqref{uniformboundsol} and \eqref{expl2} together with Lemma \ref{exp}, we obtain \eqref{decayhhalf}.

\section{Subcritical Periodic Case} \la{sub}

In this section, we consider the subcritical case where the dissipation is given by $\Lambda^{\alpha}$   
for $\alpha \in (1,2]$, that is, we consider the equation 
\be 
\pa_t \rho + u \cdot \na \rho + \Lambda^{\alpha} \rho = 0
\la{pdesub}
\ee 
posed on $\mathbb{T}^2$, where
\be 
\la{defusub}
u = - \PP (\rho R \rho).
\ee 
The initial data are given by
\be
\la{initialsub} 
\rho(x,0) = \rho_0(x)
\ee
and have zero mean. 

Global weak solutions exist:

\begin{thm} \la{Weaksub}
Let $\alpha \in (1,2]$. Let $T > 0$ be arbitrary. Let $\rho_0 \in L^2(\mathbb{T}^2)$. Then \eqref{pdesub}--\eqref{initialsub} has a weak solution $\rho$ on $[0,T]$ obeying 
\be  
\fr{1}{2} \|\rho (t)\|_{L^2 (\mathbb{T}^2)}^2 + \int_{0}^{t} \|\Lambda^{\fr{\alpha}{2}} \rho (s)\|_{L^2(\mathbb{T}^2)}^2 ds   \leq \fr{1}{2} \|\rho_0\|_{L^2(\mathbb{T}^2)}^2
\la{regsub}
\ee 
for $t \in [0,T]$.
\end{thm}
The proof is similar to that of Theorem \ref{Weak1}, and so we omit the details. 

We note that the regularity of the initial data imposed in the critical case $(\alpha = 1)$, namely $\rho_0 \in L^{2+\delta}$ for some $\delta > 0$, is not required in the subcritical case in view of the fact that $\rho$  obeys  
\be 
\rho \in L^{2}(0,T; H^{\fr{\alpha}{2}}(\mathbb{T}^2)).
\ee 

The following proposition is the analogue of Proposition \ref{crit}:

\begin{prop} \la{subcrit}
Let $\alpha \in (1,2]$. Let $p > 2$ and $\rho_0 \in L^p (\mathbb{T}^2)$. Suppose $\rho$ is a smooth solution of \eqref{pdesub}--\eqref{initialsub} on $[0,T]$. Then 
\be 
\|\rho(t)\|_{L^p (\mathbb{T}^2)} \leq \|\rho_0\|_{L^p(\mathbb{T}^2)}
\ee 
holds for all $t \in [0,T]$. Moreover, if $\rho_0 \in L^{\infty} (\mathbb{T}^2)$, then 
\be 
\|\rho(t)\|_{L^{\infty}(\mathbb{T}^2)} \leq \|\rho_0\|_{L^{\infty}(\mathbb{T}^2)} 
\ee 
holds for all $t \in [0,T]$.
\end{prop}

The solution of the initial value problem \eqref{pdesub}--\eqref{initialsub} with large smooth data are globally regular.
\beg{thm}
Let $\alpha \in (1,2]$, $s > 0$. Let $T > 0$ be arbitrary. Let $\rho_0 \in H^{s} (\mathbb{T}^2) \cap L^{\infty} (\mathbb{T}^2)$. Then  there are positive constants $C_1$, $C_2$ and $C_3$ depending only on $\|\rho_0\|_{L^{\infty}  (\mathbb{T}^2)}$ such that the solution of \eqref{pdesub}--\eqref{defusub} with initial data $\rho_0$ exists and satisfies 
\be \la{highsub1}
\|\l^{s} \rho (t)\|_{L^2(\mathbb{T}^2)} \le \|\l^{s} \rho_0\|_{L^2(\mathbb{T}^2)}e^{C_1t}
\ee and 
\be \la{highsub2}
\int_{0}^{t} \|\l^{s+\fr{\alpha}{2}} \rho (\tau) \|_{L^2(\mathbb{T}^2)}^2 d\tau
\le \|\l^{s} \rho_0\|_{L^2(\mathbb{T}^2)}^2 + C_2 \|\l^{s} \rho_0\|_{L^2(\mathbb{T}^2)}^2  (e^{C_3t} -1)
\ee
for $t \in [0,T]$.
\end{thm}

\textbf{Proof:} Fix a small $\epsilon \in (0,1)$ such that $\alpha \ge \epsilon +1$.
We multiply \eqref{pdesub} by $\l^{2s} \rho$ and we integrate in the space variable over $\mathbb{T}^2$. We obtain the equation
\be 
\fr{1}{2} \fr{d}{dt} \|\l^{s} \rho\|_{L^2(\mathbb{T}^2)}^2 + \|\l^{s + \fr{\alpha}{2}} \rho\|_{L^2(\mathbb{T}^2)}^2 
= - \int_{\mathbb{T}^2} (u \cdot \na \rho) \l^{2s} \rho dx.
\ee

We estimate the nonlinear term. Integrating by parts and using H\"older's inequality, we have 
\beg{align} 
\left|\int_{\mathbb{T}^2} (u \cdot \na \rho) \l^{2s} \rho dx \right|
&= \left|\int_{\mathbb{T}^2} \l^{s - \fr{\alpha}{2}} \na \cdot (u \rho) \l^{s + \fr{\alpha}{2}} \rho dx \right| \nonumber
\\&\le \|\l^{s- \fr{\alpha}{2} + 1} (u\rho)\|_{L^2(\mathbb{T}^2)} \|\l^{s + \fr{\alpha}{2}}\rho\|_{L^2(\mathbb{T}^2)}.
\end{align}

In view of the commutator estimate 
\be \la{commutator}
\|\l^{s} (fg)\|_{L^p(\mathbb{T}^2)} \le C\|g\|_{L^{p_1}(\mathbb{T}^2)} \|\l^s f\|_{L^{p_2}(\mathbb{T}^2)} 
+ C\|\l^s g\|_{L^{p_3}(\mathbb{T}^2)} \|f\|_{L^{p_4}(\mathbb{T}^2)}
\ee that holds for any mean zero functions $f, g \in C^{\infty}(\mathbb{T}^2), s >0, p \in (1, \infty)$ with $\fr{1}{p} = \fr{1}{p_1} + \fr{1}{p_2} = \fr{1}{p_3} + \fr{1}{p_4}, p_2, p_3 \in (1, \infty)$ (see \cite{CTV}), we estimate 
\be 
\|\l^{s- \fr{\alpha}{2} + 1} (u\rho)\|_{L^2(\mathbb{T}^2)}  
\le C\|u\|_{L^{\fr{2}{\epsilon}} (\mathbb{T}^2)} \|\l^{s - \fr{\alpha}{2} +1 }\rho\|_{L^{\fr{2}{1-\epsilon}}(\mathbb{T}^2)}
+ C\|\rho\|_{L^{\infty}(\mathbb{T}^2)} \|\l^{s-\fr{\alpha}{2} + 1}u \|_{L^2(\mathbb{T}^2)}.
\ee

In view of the boundedness of the Riesz transforms (and hence the Leray projector) on $L^p (\mathbb{T}^2)$ for $p \in (1, \infty)$ and Proposition \ref{subcrit}, we bound
\be 
\|u\|_{L^{\fr{2}{\epsilon}}(\mathbb{T}^2)} 
\le C\|\rho R \rho\|_{L^{\fr{2}{\epsilon}}(\mathbb{T}^2)} 
\le C\|\rho\|_{L^{\infty} (\mathbb{T}^2)} \|\rho\|_{L^{\fr{2}{\epsilon}}(\mathbb{T}^2)}
\le C \|\rho\|_{L^{\infty}(\mathbb{T}^2)}^2
\le C\|\rho_0\|_{L^{\infty}(\mathbb{T}^2)}^2.
\ee

By the commutator estimate \eqref{commutator}, we have
\beg{align}
&\|\l^{s - \fr{\alpha}{2} + 1} u\|_{L^2(\mathbb{T}^2)}
\le C\|\l^{s - \fr{\alpha}{2} + 1} (\rho R \rho)\|_{L^2(\mathbb{T}^2)} \nonumber
\\&\le C\|\rho\|_{L^{\infty}(\mathbb{T}^2)} \|\l^{s - \fr{\alpha}{2}+1} R\rho \|_{L^2 (\mathbb{T}^2)}
+ C\|R\rho\|_{L^{\fr{2}{\epsilon}}(\mathbb{T}^2)} \|\l^{s - \fr{\alpha}{2} + 1} \rho \|_{L^\fr{2}{1-\epsilon}(\mathbb{T}^2)} \nonumber
\\&\le C\|\rho_0\|_{L^{\infty}(\mathbb{T}^2)} \|\l^{s - \fr{\alpha}{2}+1} \rho \|_{L^2(\mathbb{T}^2)}
+ C\|\rho_0\|_{L^{\infty}(\mathbb{T}^2)} \|\l^{s - \fr{\alpha}{2} + 1} \rho \|_{L^\fr{2}{1-\epsilon}(\mathbb{T}^2)} 
\end{align}

Hence
\be 
\|\l^{s- \fr{\alpha}{2} + 1} (u\rho)\|_{L^2(\mathbb{T}^2)}  
\le C\|\rho_0\|_{L^{\infty}(\mathbb{T}^2)}^2 \|\l^{s - \fr{\alpha}{2} +1 }\rho\|_{L^{\fr{2}{1-\epsilon}}(\mathbb{T}^2)}
+ C\|\rho_0\|_{L^{\infty}(\mathbb{T}^2)}^2 \|\l^{s - \fr{\alpha}{2}+1} \rho \|_{L^2(\mathbb{T}^2)}.
\ee

In view of the continuous Sobolev embedding 
\be 
H^{\epsilon} (\mathbb{T}^2) \hookrightarrow L^{\fr{2}{1-\epsilon}} (\mathbb{T}^2),
\ee we obtain the bound
\be 
\|\l^{s- \fr{\alpha}{2} + 1} (u\rho)\|_{L^2(\mathbb{T}^2)}  
\le C\|\rho_0\|_{L^{\infty}(\mathbb{T}^2)}^2 \|\l^{s - \fr{\alpha}{2} +1 + \epsilon}\rho\|_{L^{2}(\mathbb{T}^2)}
+ C\|\rho_0\|_{L^{\infty}(\mathbb{T}^2)}^2 \|\l^{s - \fr{\alpha}{2}+1} \rho \|_{L^2(\mathbb{T}^2)}.
\ee

Using the Sobolev interpolation inequality 
\be 
\|\l^{s_1} f\|_{L^2(\mathbb{T}^2)} \le C\|\l^{s_0} f\|_{L^2(\mathbb{T}^2)}^{1-\sigma} \|\l^{s_2} f\|_{L^2(\mathbb{T}^2)}^{\sigma}
\ee that holds for any mean zero function $f \in H^{s_2}(\mathbb{T}^2)$ and $s_1 = (1-\sigma)s_0 + \sigma s_2$, $\sigma \in [0,1]$, we estimate
\be 
\|\l^{s - \fr{\alpha}{2}+1} \rho \|_{L^2(\mathbb{T}^2)}
\le C\left(\|\l^{s} \rho\|_{L^2(\mathbb{T}^2)}\right)^{\fr{2(\alpha -1)}{\alpha}} \left(\|\l^{s+\fr{\alpha}{2}} \rho\|_{L^2(\mathbb{T}^2)} \right)^{\fr{2}{\alpha} -1}
\ee 
and 
\be 
\|\l^{s - \fr{\alpha}{2} +1 + \epsilon}\rho\|_{L^{2}(\mathbb{T}^2)}
\le C \left(\|\l^{s} \rho\|_{L^2(\mathbb{T}^2)}\right)^{\fr{2(\alpha -\epsilon -1)}{\alpha}} \left(\|\l^{s+\fr{\alpha}{2}} \rho\|_{L^2(\mathbb{T}^2)} \right)^{\fr{2(\epsilon +1)}{\alpha}-1}
\ee

Consequently, 
\beg{align} 
&\|\l^{s- \fr{\alpha}{2} + 1} (u\rho)\|_{L^2(\mathbb{T}^2)} \|\l^{s + \fr{\alpha}{2}}\rho\|_{L^2(\mathbb{T}^2)}  \nonumber
\\&\le C\|\rho_0\|_{L^{\infty}(\mathbb{T}^2)}^2 \left(\|\l^{s} \rho\|_{L^2(\mathbb{T}^2)}\right)^{\fr{2(\alpha -\epsilon -1)}{\alpha}} \left(\|\l^{s+\fr{\alpha}{2}} \rho\|_{L^2(\mathbb{T}^2)} \right)^{\fr{2(\epsilon +1)}{\alpha}} \nonumber
\\&+ C\|\rho_0\|_{L^{\infty}(\mathbb{T}^2)}^2 \left(\|\l^{s} \rho\|_{L^2(\mathbb{T}^2)}\right)^{\fr{2(\alpha -1)}{\alpha}} \left(\|\l^{s+\fr{\alpha}{2}} \rho\|_{L^2(\mathbb{T}^2)} \right)^{\fr{2}{\alpha}}
\end{align}

By Young's inequality, we end up with
\be 
\left|\int_{\mathbb{T}^2} (u \cdot \na \rho) \l^{2s} \rho dx \right|
\le C_{\rho_0} \|\l^{s} \rho\|_{L^2(\mathbb{T}^2)}^2 + \fr{1}{2} \|\l^{s+\fr{\alpha}{2}} \rho\|_{L^2(\mathbb{T}^2)}^2
\ee
where $C_{\rho_0}$ is a constant depending on the $L^{\infty}$ norm of the initial data $\rho_0$. 

Therefore, we obtain the differential inequality,
\be 
\fr{d}{dt} \|\l^{s} \rho\|_{L^2(\mathbb{T}^2)}^2 + \|\l^{s + \fr{\alpha}{2}} \rho\|_{L^2(\mathbb{T}^2)}^2 
\le 2C_{\rho_0} \|\l^{s} \rho\|_{L^2(\mathbb{T}^2)}^2
\ee
which gives \eqref{highsub1} and \eqref{highsub2}. 

We have shown existence of global smooth solutions in the subcritical case, provided that the initial data is smooth enough. No smallness condition is imposed on the size of the initial data. The solutions are also unique. The results obtained hold as well in the whole space $\R^2$.

\section{Appendix A} \la{prop2}

In this appendix, we prove Proposition \ref{Bon}.  Let $f, g \in \mathcal{S}_h'$. Bony's paraproduct gives the decomposition  
\be 
fg 
= \sum\limits_{j \in \ZZ} S_{j-1}f \Delta_j g 
+ \sum\limits_{j \in \ZZ} S_{j-1}g \Delta_j f 
+ \sum\limits_{|j - j'| \le 1} \Delta_j f \Delta_{j'}g.
\ee 
We note that 
\beg{align} 
\sum\limits_{|j - j'| \le 1} \Delta_j f \Delta_{j'}g
&= \sum\limits_{j \in \ZZ} \Delta_j f \Delta_j g 
+ \sum\limits_{j \in \ZZ} \Delta_j f \Delta_{j-1} g
+ \sum\limits_{j \in \ZZ} \Delta_j f \Delta_{j+1} g \nonumber 
\\&=  \sum\limits_{j \in \ZZ} \Delta_j f \Delta_j g 
+ \sum\limits_{j \in \ZZ} \Delta_j f \Delta_{j-1} g
+ \sum\limits_{j \in \ZZ} \Delta_{j-1} f \Delta_{j} g \nonumber
\\&=  \sum\limits_{j \in \ZZ} (\Delta_{j-1} f  + \Delta_j f) \Delta_j g 
+ \sum\limits_{j \in \ZZ} \Delta_j f \Delta_{j-1} g.  
\end{align}
This implies that
\be 
fg = \sum\limits_{j \in \ZZ} S_{j+1}f \Delta_j g 
+  \sum\limits_{j \in \ZZ} S_{j}g \Delta_j f 
\ee
Now we apply $\Delta_j$. In view of \eqref{supports}, we have
\be \la{supp1}
k \le j - 2 \Rightarrow \Delta_j (S_{k} g \Delta_k f) = 0
\ee and 
\be  \la{supp2}
k \le j - 3 \Rightarrow \Delta_j (S_{k+1} f \Delta_k g) = 0
\ee  
Indeed, \beg{align}
&\mathcal{F} (\Delta_j (S_k g\Delta_k f) (\xi)
= \Psi_j(|\xi|) \mathcal{F} (S_kg \Delta_k f)(\xi) \nonumber
\\&=  \Psi_j (|\xi|) \left\{\sum\limits_{l \le k-1} \int\limits_{\R^2} \Psi_{l}(|\xi - y|) \mathcal{F}g(\xi - y) \Psi_k (|y|) \mathcal{F}f(y) dy \right\} \nonumber
\\&= \Psi_j (|\xi|) \left\{\sum\limits_{l \le k-1} \int\limits_{\fr{2^k}{2} \le |y| \le \fr{2^k 5}{4}} \Psi_{l}(|\xi - y|) \mathcal{F}g(\xi - y) \Psi_k (|y|) \mathcal{F}f(y) dy \right\} \nonumber
\\&= \Psi_j (|\xi|) \tilde{\Psi}_k(\xi)
\end{align} where 
\be 
\tilde{\Psi}_k(\xi) = \sum\limits_{l \le k-1} \int\limits_{\fr{2^k}{2} \le |y| \le \fr{2^k 5}{4}} \Psi_{l}(|\xi - y|) \mathcal{F}g(\xi - y) \Psi_k (|y|) \mathcal{F}f(y) dy. 
\ee Fix $l \le k-1$. Let $y \in \R^2$ such that $\fr{2^k}{2} \le |y| \le \fr{2^k 5}{4}$ and $\Psi_l (|\xi-y|) \ne 0$. This implies that $|\xi - y| \le \fr{2^l 5}{4}$, thus
\be 
|\xi|  \le |\xi - y| + |y| \le \fr{2^l5}{4} + \fr{2^k 5}{4} \le \fr{2^{k-1}5}{4} + \fr{2^k5}{4} = 2^{k-3} 15.
\ee 
Consequently, if $|\xi| > 2^{k-3}15$, then $\Psi_l (|\xi - y|) = 0$ for all $l \le k-1$ and for all $y$ satisfying $\fr{2^k}{2} \le |y| \le \fr{2^k 5}{4}$, and so $\tilde{\Psi}_k (\xi)  = 0$. We conclude that the support of $\tilde{\Psi}_k$ is included in the closed ball centered at $0$ with radius $2^{k-3}15$. But the support of $\Psi_j(|\cdot|)$ is included in the closed annulus centered at $0$ with radii $\fr{2^j}{2}$ and $\fr{2^j 5}{4}$. Therefore, if $k+1 \le j- 1$, then $2^{k-3}15 < 2^{k+1} \le 2^{j-1}$ and so 
\be 
\mathcal{F} (\Delta_j (S_k g \Delta_k f)) = 0 
\ee which gives \eqref{supp1}. The property \eqref{supp2} follows from a similar argument. 
Therefore, we obtain the decomposition  
\be 
\Delta_j (fg) = \sum\limits_{k \ge j-2} \Delta_j (S_{k+1}f \Delta_k g) 
+  \sum\limits_{k \ge j-2} \Delta_j (S_{k}g \Delta_k f). 
\ee This ends the proof of Proposition \ref{Bon}.

\section{Appendix B~\la{weakapp}}

\textbf{Proof of Theorem  \ref{Weak1}:} 
We take the $L^2$ inner product of \eqref{approx} with $\rho^{\epsilon}$  and we obtain 
\be 
\frac{1}{2} \fr{d}{dt} \|\rho^{\epsilon}\|_{L^2}^2 + \|\Lambda^{\frac{1}{2}} \rho^{\epsilon}\|_{L^2}^2 + \epsilon \|\na \rho^{\epsilon}\|_{L^2}^2 = 0.
\la{l2normeq} 
\ee 
Here we used the fact that $\widetilde{u}^{\epsilon}$ is divergence free, which implies that
\be 
(\widetilde{u}^{\epsilon} \cdot \na \rho^{\epsilon}, \rho^{\epsilon})_{L^2} = 0.
\ee 
Integrating \eqref{l2normeq} in time from $0$ to $t$, we obtain \eqref{reg1}. 
Therefore, the family $\left\{\rho^{\epsilon} : \epsilon \in (0,1] \right\}$ is uniformly bounded in $L^2(0, T; H^{\fr{1}{2}})$. 
Moreover, we have
\be 
|(\Lambda \rho^{\epsilon}, \Phi)_{L^2}| = |(\Lambda^{\fr{1}{2}}\rho^{\epsilon}, \Lambda^{\fr{1}{2}}\Phi )_{L^2}| 
\leq \|\Lambda^{\fr{1}{2}}\rho^{\epsilon} \|_{L^2} \|\Lambda^{\fr{1}{2}}\Phi \|_{L^2} 
\leq C\|\Lambda^{\fr{1}{2}}\rho^{\epsilon} \|_{L^2} \|\Phi \|_{H^{\fr{5}{2}}},
\ee 
\be 
\epsilon|(-\Delta \rho^{\epsilon}, \Phi)_{L^2}| = \epsilon |(\rho^{\epsilon}, -\Delta \Phi )_{L^2}|  
\leq C\|\rho^{\epsilon} \|_{L^2} \|\Phi \|_{H^{\fr{5}{2}}},
\ee and
\be 
|(\widetilde{u}^{\epsilon} \cdot \nabla \rho^{\epsilon}, \Phi)_{L^2}| 
= |(\widetilde{u}^{\epsilon}  \rho^{\epsilon}, \nabla  \Phi)_{L^2}|
\leq \|\widetilde{u}^{\epsilon}\|_{L^2} \|\rho^{\epsilon}\|_{L^2} \|\nabla  \Phi\|_{L^{\infty}}
\leq C\|\rho^{\epsilon}\|_{L^4}^2\|\rho^{\epsilon}\|_{L^2} \|\Phi \|_{H^{\fr{5}{2}}}
\ee 
for all $\Phi \in H^{\fr{5}{2}}$. Here we used the boundedness of the Riesz operator on $L^4$, and the continuous Sobolev embedding $H^{\fr{3}{2}} \hookrightarrow L^{\infty}$. Therefore, we obtain the bound
\be 
\|\widetilde{u}^{\epsilon} \cdot \nabla \rho^{\epsilon}\|_{H^{-\fr{5}{2}}} + \|\Lambda \rho^{\epsilon}\|_{H^{-\fr{5}{2}}} + \epsilon \|\Delta \rho^{\epsilon}\|_{H^{-\fr{5}{2}}}
\leq C(\|\rho^{\epsilon}\|_{L^4}^2\|\rho^{\epsilon}\|_{L^2} + \|\rho^{\epsilon}\|_{L^2} + \|\Lambda^{\fr{1}{2}}\rho^{\epsilon} \|_{L^2}).
\ee 
In view of the continuous embedding $H^{\fr{1}{2}} \hookrightarrow L^4$, we conclude that the family $\left\{\pa_t \rho^{\epsilon} : \epsilon \in (0,1] \right\}$ is uniformly bounded in $L^1(0,T; H^{-\fr{5}{2}})$. Now, we note that the inclusion $H^{\fr{1}{2}} \hookrightarrow L^2$ is compact whereas the inclusion $L^2 \hookrightarrow H^{-\fr{5}{2}}$ is continuous. Let $\epsilon_n$ be a decreasing sequence in $(0,1]$ converging to 0. By the Aubin-Lions lemma and \eqref{reg1}, the sequence $\left\{\rho^{\epsilon_n}\right\}_{n=1}^{\infty}$ has a subsequence  that converges strongly in $L^2(0,T; L^2)$ and weakly in $L^2(0,T; H^{\fr{1}{2}})$ to some function $\rho$. By the lower semi-continuity of the norms, we obtain \eqref{reg2}. 

For simplicity of notations, we assume that $\rho^{\epsilon}$ converges to $\rho$ strongly in $L^2(0,T; L^2)$ and weakly in $L^2(0,T; H^{\fr{1}{2}})$. 
We note that 
\be 
(\rho^{\epsilon} (t), \Phi)_{L^2} - (\rho_0, \Phi)_{L^2} 
+ \int_{0}^{t} (\widetilde{u}^{\epsilon} \cdot \nabla \rho^{\epsilon}, \Phi)_{L^2}  ds 
+ \int_{0}^{t} (\Lambda^{\fr{1}{2}} \rho^{\epsilon}, \Lambda^{\fr{1}{2}}  \Phi)_{L^2} ds 
+ \epsilon \int_{0}^{t} (\na \rho^{\epsilon}, \na \Phi)_{L^2} ds
= 0
\ee holds for all $\Phi \in H^{\fr{5}{2}}$ and $t \in [0,T]$.
Without loss of generality, we may assume that $\rho^{\epsilon}$ converges $\rho$ in $L^2$ for a.e. $t \in [0,T]$, and so 
\be 
|(\rho^{\epsilon} (t), \Phi)_{L^2} - (\rho(t), \Phi)_{L^2}| \leq \|\rho^{\epsilon} - \rho\|_{L^2} \|\Phi\|_{L^2} \rightarrow 0
\ee for all $\Phi \in H^{\fr{5}{2}}$ and a.e. $t \in [0,T]$. By the weak convergence in $L^2(0,T; H^{\fr{1}{2}})$, we obtain
\be 
\left|\int_{0}^{t} (\Lambda^{\fr{1}{2}} \rho^{\epsilon}, \Lambda^{\fr{1}{2}}  \Phi)_{L^2} ds - \int_{0}^{t} (\Lambda^{\fr{1}{2}} \rho, \Lambda^{\fr{1}{2}}  \Phi)_{L^2} ds \right| \rightarrow 0
\ee for all $\Phi \in H^{\fr{5}{2}}$ and all $t \in [0,T]$. 
For the nonlinear term, we let $\Phi \in H^{\fr{5}{2}}$, $t \in [0,T]$ and we write 
\beg{align}
&\int_{0}^{t} (\widetilde{u}^{\epsilon} \cdot \nabla \rho^{\epsilon}, \Phi)_{L^2}  ds - \int_{0}^{t} (u \cdot \nabla \rho , \Phi)_{L^2}  ds  \nonumber
\\&= -\int_{0}^{t} ((\rho^{\epsilon} - \rho)u, \nabla \Phi)_{L^2} ds 
- \int_{0}^{t} ((\widetilde{u}^{\epsilon} - u) \rho^{\epsilon}, \nabla \Phi)_{L^2}  ds  \nonumber
\\&= I_1 + I_2.
\end{align} We note that 
\be 
|I_1| \leq C\|\Phi\|_{H^{\fr{5}{2}}} \int_{0}^{t} \|\rho\|_{L^4}^2 \|\rho^{\epsilon} - \rho\|_{L^2} ds \rightarrow 0
\ee by the Lebesgue Dominated Convergence theorem. For $I_2$, we split it as 
\beg{align} 
I_2 &=  \int_{0}^{t} ((J_{\epsilon} \PP(\rho (R \rho^{\epsilon} - R \rho))) \rho^{\epsilon}, \nabla \Phi)_{L^2}  ds  
+ \int_{0}^{t} ((J_{\epsilon} \PP ((\rho^{\epsilon} - \rho)R\rho^{\epsilon})) \rho^{\epsilon}, \nabla \Phi)_{L^2}  ds \nonumber
\\&= I_{2,1} + I_{2,2}.  
\end{align}
In view of the boundedness of the Riesz transform on $L^2$ and the boundedness of the Leray operator on $L^{4/3}$, we have \beg{align}
|I_{2,1}| &\leq C\|\Phi\|_{H^{\fr{5}{2}}} \int_{0}^{t} \|\rho^{\epsilon}\|_{L^4} \|\PP (\rho R(\rho^{\epsilon} - \rho))\|_{L^{4/3}}  ds \nonumber
\\& \leq C\|\Phi\|_{H^{\fr{5}{2}}} \int_{0}^{t} \|\rho^{\epsilon}\|_{L^4} \|\rho\|_{L^4} \|\rho^{\epsilon} - \rho\|_{L^2} ds   \nonumber
\\&\leq C\|\Phi\|_{H^{\fr{5}{2}}} \left(\int_{0}^{t} \|\rho^{\epsilon}\|_{L^4}^2  ds \right)^{1/2} \left(\int_{0}^{t} \|\rho\|_{L^4}^2 \|\rho^{\epsilon} - \rho\|_{L^2}^2 ds \right)^{1/2} 
\rightarrow 0
\end{align} by the Lebesgue Dominated Convergence theorem. 

We note that we have not yet used the assumption that $\rho_0 \in L^{2+ \delta}$. It will be needed to estimate $|I_{2,2}|$. Indeed, we multiply equation \eqref{approx} by $\rho^{\epsilon}|\rho^{\epsilon}|^{\delta}$ and we integrate in the space variable. We use the C\'ordoba-C\'ordoba inequality \cite{CC}  
\beg{align}
&\int_{\mathbb{T}^2} |\rho^{\epsilon}|^{\delta} (\rho^{\epsilon} \Lambda{\rho^{\epsilon}}) dx  
\geq 0
\end{align}
and we obtain the differential inequality  
\be 
\fr{d}{dt} \|\rho^{\epsilon}(t)\|_{L^{2 + \delta}} \le 0.
\ee
Integrating in time from $0$ to $t$, we end up having the bound 
\be 
\|\rho^{\epsilon}(t)\|_{L^{2+\delta}} \leq \|\rho_0\|_{L^{2+\delta}}
\ee 
for all $t \in [0,T]$.
As a consequence, 
\beg{align} 
|I_{2,2}| &\leq C\|\Phi\|_{H^{\fr{5}{2}}} \int_{0}^{t}  \|\rho^{\epsilon}\|_{L^{4}} \|\rho^{\epsilon}\|_{L^{2+\delta}}  \|\rho^{\epsilon} - \rho\|_{L^{\fr{8+4\delta}{2 + 3\delta}}} ds  \nonumber
\\&\leq C\|\Phi\|_{H^{\fr{5}{2}}} \|\rho_0\|_{L^{2+\delta}} \int_{0}^{t} \|\rho^{\epsilon}\|_{L^4} \|\rho^{\epsilon} - \rho\|_{L^2}^{\fr{2\delta}{2 + \delta}}  \|\rho^{\epsilon} - \rho\|_{L^4}^{\fr{2 - \delta}{2 + \delta}}  ds \nonumber
\\&\leq C\|\Phi\|_{H^{\fr{5}{2}}} \|\rho_0\|_{L^{2+\delta}} \left(\int_{0}^{t} \|\rho^{\epsilon}\|_{L^4}^2 \right)^{\fr{2}{2 + \delta}}\left(\int_{0}^{t} \|\rho^{\epsilon} - \rho\|_{L^2}^{2} ds \right)^{\fr{\delta}{2+\delta}} \nonumber
\\&+ C\|\Phi\|_{H^{\fr{5}{2}}} \|\rho_0\|_{L^{2+\delta}} \left(\int_{0}^{t} \|\rho^{\epsilon}\|_{L^4}^2 \right)^{1/2} \left(\int_{0}^{t} \|\rho\|_{L^4}^{2} ds \right)^{\fr{2 - \delta}{4 + 2\delta}} \left(\int_{0}^{t} \|\rho^{\epsilon} - \rho\|_{L^2}^{2} ds \right)^{\fr{\delta}{2+\delta}}  
\rightarrow 0.
\end{align}
Here we used the interpolation inequality 
\be 
\|f\|_{L^{\fr{8 + 4\delta}{2 + 3 \delta}}} \leq C\|f\|_{L^2}^{\fr{2\delta}{2 + \delta}} \|f\|_{L^4}^{\fr{2 - \delta}{2 + \delta}}
\ee 
that holds for any $f \in L^{4}$. 

Therefore $\rho$ is a weak solution of \eqref{pde}. This ends the proof of Theorem \ref{Weak1}.

\section{Appendix C \la{appstrong} }
\textbf{Proof of Theorem \ref{strongloc}:} We apply $-\Delta = \Lambda^2$ to \eqref{approx} and we obtain
\be 
-\partial_t \Delta \rho^{\epsilon} - \widetilde{u}^{\epsilon} \cdot \nabla \Delta \rho^{\epsilon} -2\nabla \widetilde{u}^{\epsilon}  \nabla \nabla \rho^{\epsilon} - \Delta \widetilde{u}^{\epsilon} \cdot \nabla \rho^{\epsilon} + \Lambda^3 \rho^{\epsilon} + \epsilon \Delta \Delta \rho^{\epsilon} = 0 
\la{deltar}
\ee 
We multiply \eqref{deltar} by $-\Delta \rho^{\epsilon}$ and we integrate over $\mathbb{R}^2$. In view of the fact that
\be 
 (\widetilde{u}^{\epsilon} \cdot \nabla \Delta \rho^{\epsilon}, \Delta \rho^{\epsilon})_{L^2} = 0,
\ee 
we obtain
\be 
\fr{1}{2} \fr{d}{dt} \|\Delta \rho^{\epsilon}\|_{L^2}^2 + \|\Lambda^{\fr{5}{2}} \rho^{\epsilon}\|_{L^2}^2  
+ \epsilon \|\l^3 \rho^{\epsilon}\|_{L^2}^2
= - 2(\nabla \widetilde{u}^{\epsilon} \nabla \nabla \rho^{\epsilon} ,\Delta \rho^{\epsilon})_{L^2} - (\Delta \widetilde{u}^{\epsilon} \cdot \nabla \rho^{\epsilon},\Delta \rho^{\epsilon})_{L^2}.
\ee  

Using the product rule 
\be 
\|fg\|_{H^s} \le C\|f\|_{H^s}\|g\|_{L^{\infty}} + C\|g\|_{H^s}\|f\|_{L^{\infty}}
\ee that holds for any $f, g \in H^s, s > 0$, we estimate 
\beg{align}
&\|\nabla \widetilde{u}^{\epsilon}\|_{L^4} 
\le C \| \widetilde{u}^{\epsilon}\|_{H^{\fr{3}{2}}} 
\le C\| \rho^{\epsilon} R \rho^{\epsilon} \|_{H^{\fr{3}{2}}} \nonumber
\\&\le C \|\rho^{\epsilon}\|_{L^{\infty}} \|R \rho^{\epsilon}\|_{H^{\fr{3}{2}}}
+ C \|R\rho^{\epsilon}\|_{L^{\infty}}\| \rho^{\epsilon}\|_{H^{\fr{3}{2}}} \nonumber
\\&\le C \|\rho^{\epsilon}\|_{H^{\fr{3}{2}}}^2.
\end{align} Here, we have used the continuous embedding $H^{\fr{1}{2}}\hookrightarrow L^4$, the fact that the Leray projector is bounded on $H^{\fr{3}{2}}$, and the boundedness of the Riesz transforms as operators from $H^{\fr{3}{2}}$  into $L^{\infty}$. 
Similarly, we bound 
\beg{align}
&\|\Delta \widetilde{u}^{\epsilon}\|_{L^4} 
\le C\|\rho^{\epsilon} R \rho^{\epsilon}\|_{H^{\fr{5}{2}}}  \nonumber
\\&\le C \|\rho^{\epsilon}\|_{L^{\infty}} \|R \rho^{\epsilon}\|_{H^{\fr{5}{2}}}
+ C \|R\rho^{\epsilon}\|_{L^{\infty}}\|R \rho^{\epsilon}\|_{H^{\fr{5}{2}}} \nonumber
\\&\le C \|\rho^{\epsilon}\|_{H^{\fr{3}{2}}} \|\rho^{\epsilon}\|_{H^{\fr{5}{2}}}
\end{align}
Consequently, 
\beg{align}
\fr{1}{2} \fr{d}{dt} \|\Delta \rho^{\epsilon}\|_{L^2}^2 + \|\Lambda^{\fr{5}{2}} \rho^{\epsilon}\|_{L^2}^2  
&\le 2 \|\nabla \widetilde{u}^{\epsilon}\|_{L^4} \|\na \na \rho^{\epsilon}\|_{L^4} \|\Delta \rho^{\epsilon}\|_{L^2}
+ \|\Delta \widetilde{u}^{\epsilon}\|_{L^4}  \|\na \rho^{\epsilon}\|_{L^4} \|\Delta \rho^{\epsilon}\|_{L^2} \nonumber
\\&\le C\|\rho^{\epsilon}\|_{H^{\fr{3}{2}}}^2 \|\rho^{\epsilon}\|_{H^{\fr{5}{2}}} \|\Delta \rho^{\epsilon}\|_{L^2}
\end{align}
and by Young's inequality, we obtain 
\beg{align} 
&\fr{d}{dt} \|\Delta \rho^{\epsilon}\|_{L^2}^2 + \|\Lambda^{\fr{5}{2}} \rho^{\epsilon}\|_{L^2}^2 
\le C\|\rho^{\epsilon}\|_{H^{\fr{3}{2}}}^4 \|\Delta \rho^{\epsilon}\|_{L^2}^2 
+ C\|\rho^{\epsilon}\|_{H^{\fr{3}{2}}}^2 \|\rho^{\epsilon}\|_{L^2} \|\Delta \rho^{\epsilon}\|_{L^2} \nonumber
\\&\le C(\|\rho^{\epsilon}\|_{H^2}^6 + \|\rho^{\epsilon}\|_{H^2}^4).
\end{align}
We note that 
\beg{align} 
&\|\rho^{\epsilon}\|_{H^2} = \left\|(1+|.|^{2})\mathcal{F}(\rho^{\epsilon})(.) \right\|_{L^2}
\le C\|\mathcal{F} \rho^{\epsilon}\|_{L^2} + C \|\Delta \rho^{\epsilon}\|_{L^2} \nonumber
\\&= C\|\rho^{\epsilon}\|_{L^2} + C \|\Delta \rho^{\epsilon}\|_{L^2}
\le C\|\rho_0\|_{L^2} + C\|\Delta \rho^{\epsilon}\|_{L^2}
\end{align} in view of Plancherel's theorem and the uniform boundedness of $\rho^{\epsilon}$ in $L^2$ described by \eqref{reg2}. Therefore, we obtain the differential inequality
\be 
\fr{d}{dt} \|\Delta \rho^{\epsilon}\|_{L^2}^2 + \|\Lambda^{\fr{5}{2}} \rho^{\epsilon}\|_{L^2}^2 
\le C \|\Delta \rho^{\epsilon}\|_{L^2}^6 + C_{\rho_0} 
\ee where $C_{\rho_0}$ is a positive constant depending only on $\rho_0$ and some universal constants. This gives a local strong solution. 

For uniqueness, suppose that $\rho_1$ and $\rho_2$ are two strong solutions of \eqref{pde} on $[0,T_0]$ with the same initial condition. Let $\rho = \rho_1 - \rho_2$ and $u = u_1 - u_2$. Then $\rho$ obeys the equation 
\be 
\partial_t \rho + u \cdot \nabla \rho_1 + u_2 \cdot \nabla \rho + \Lambda \rho = 0
\ee  
We take the $L^2$ inner product with $\rho$ and we obtain 
\be 
\fr{1}{2} \fr{d}{dt} \|\rho\|_{L^2}^2 + \|\Lambda^{\fr{1}{2}}\rho \|_{L^2}^2 = -(u \cdot \nabla \rho_1, \rho)_{L^2}.
\ee

In view of the boundedness of the Riesz transforms on $L^4$, we have
\beg{align}
\|u\|_{L^4} &\le \|\mathbb{P}(\rho R\rho_1)\|_{L^4} + \|\mathbb{P}(\rho_2 R\rho)\|_{L^4} \nonumber
\\&\le C\|\rho\|_{L^{4}} \|R\rho_1\|_{L^{\infty}} + \|\rho_2\|_{L^{\infty}}\|R\rho\|_{L^{4}} \nonumber
\\&\le C\|\rho\|_{L^4} \left(\|\rho_1\|_{H^{\fr{3}{2}}} + \|\rho_2\|_{H^{\fr{3}{2}}} \right).
\end{align}
Hence
\beg{align}
&|(u \cdot \nabla \rho_1, \rho)_{L^2} | 
\le \|u\|_{L^4} \|\nabla \rho_1\|_{L^4} \|\rho\|_{L^2}  \nonumber
\\&\le \fr{1}{2} \|\rho\|_{H^{\fr{1}{2}}}^2 
+ C\left(\|\rho_1\|_{H^{\fr{3}{2}}}^2 + \|\rho_2\|_{H^{\fr{3}{2}}}^2\right) \|\rho_1\|_{H^{\fr{3}{2}}}^2\|\rho\|_{L^2}^2.
\end{align}
Therefore, 
\be 
\fr{d}{dt}\|\rho\|_{L^2}^2 \leq K(t) \|\rho\|_{L^2}^2
\ee
where 
\be 
K(t) = C\left(\|\rho_1\|_{H^{\fr{3}{2}}}^2 + \|\rho_2\|_{H^{\fr{3}{2}}}^2\right) \|\rho_1\|_{H^{\fr{3}{2}}}^2.
\ee
This shows that for each $t \geq 0$, $\rho_1(\cdot,t) = \rho_2(\cdot,t)$ a.e. in $\mathbb{R}^2$.

\end{document}